\documentclass[10pt]{amsart}
\usepackage{amssymb, graphicx, labelfig, subfig
}

\newtheorem{thm}{Theorem}
\newtheorem{prop}[thm]{Proposition}
\newtheorem{lem}[thm]{Lemma}

\theoremstyle{remark}
\newtheorem{rem}[thm]{Remark}

\theoremstyle{definition}

\newcommand{\col}{\kern -3pt :}
\newcommand{\C}{\mathbb C}
\newcommand{\R}{\mathbb R}
\newcommand{\Z}{\mathbb Z}

\newcommand{\HH}{\mathbb H}

\newcommand{\Id}{\mathrm{Id}}

\newcommand{\Spin}{\mathrm{Spin}}
\newcommand{\RS}{\mathcal R_{\PSL(\C)}^{\Spin}}

\newcommand{\T}{\mathcal T}
\newcommand{\SSS}{\mathcal S^A}
\newcommand{\SSSS}{\mathcal S_{\mathrm s}^A}
\newcommand{\RR}{\mathcal R}
\newcommand{\ZZ}{\mathcal Z^\omega}
\newcommand{\ZZZ}{\widehat{\mathcal Z}^\omega}
\newcommand{\SL}{\mathrm{SL}_2}
\newcommand{\PSL}{\mathrm{PSL}_2}
\newcommand{\Tr}{\mathrm{Tr}}
\newcommand{\E}{\mathrm{e}}

\newcommand{\db}{/\kern -4pt/}

\renewcommand{\geq}{\geqslant}
\renewcommand{\phi}{\varphi}
\renewcommand{\epsilon}{\varepsilon}

\title[Quantum traces]
{Quantum traces for representations of surface groups in $\SL$}

\author{Francis Bonahon}
\address {Department
of Mathematics,  University of
Southern California, Los Angeles
CA~90089-2532, U.S.A.}
\email{fbonahon@math.usc.edu}

\author{Helen Wong}
\address{Department
of Mathematics, Carleton College, Northfield MN 55057, U.S.A.}
\email{hwong@carleton.edu}

\thanks{This research was partially supported by grants DMS-0632713 (Carleton Summer Mathematics Program) and DMS-0604866 from the National Science Foundation.}
\date{\today}

\begin{document}

\setlength{\baselineskip 13pt}
\maketitle

\begin{abstract}
We consider two different quantizations of the character variety consisting of all representations of surface groups in $\PSL$. One is the skein algebra considered by Bullock-Frohman-Kania-Bartoszy\'nska, Przytycki-Sikora and Turaev. The other is the quantum Teichm\"uller space introduced by Chekhov-Fock and Kashaev.  We construct a homomorphism from the skein algebra to the quantum Teichm\"uller space which, when restricted the classical case, corresponds to the equivalence between these two algebras through trace functions. 
\end{abstract}

Let $S$ be an oriented surface of finite topological type. The goal of this paper is to establish a connection between two quantizations of the character  variety
$$\RR_{\SL(\C)}(S) = \{ r \col \pi_1(S) \to \SL(\C)\}\db\SL(\C),$$
consisting of all group homomorphisms $r$ from the fundamental group $\pi_1(S)$ to the Lie group $\SL(\C)$, considered up to conjugation by elements of $\SL(\C)$. The double bar indicates here that the quotient is taken in the algebraic geometric sense of geometric invariant theory.

The first quantization, introduced by  D.~Bullock, C.~Frohman, J.~Kania-Bartos\-zy\'nska \cite{BFK1}, J.~Przytycki, A.~Sikora \cite{PrzS} and V.~Turaev \cite{Tur}, uses the \emph{skein algebra} $\SSS(S)$ obtained by considering the vector space freely generated by all isotopy classes of framed links in $S\times [0,1]$, and then taking the quotient of this space under the Kauffman skein relation; see \S \ref{sect:skein}. What makes $\SSS(S)$ a quantization of $\RR_{\SL(\C)}(S)$ is that, when $A=-1$, the skein algebra $\mathcal S^{-1}(S)$ has a natural identification with the commutative algebra of regular functions on $\RR_{\SL(\C)}(S)$ and  that, as $A$ tends to $-1$, the lack of commutativity of $\SSS(S)$ is infinitesimally measured by the Goldman-Weil-Petersson Poisson structure \cite{Gold1, Gold, Pet, Weil} on $\RR_{\SL(\C)}(S)$; see \cite{Tur}. There is a similar situation when $A=+1$, in which case $\mathcal S^{+1}(S)$ has a natural identification with a twisted version of $\RR_{\SL(\C)}(S)$; see \S \ref{subsect:ClassicalSkeins}. 

The second quantization, with respect to the same Goldman-Weil-Petersson Poisson structure, is the quantum Teichm\"uller space $\widehat \T_S^q$ introduced by V.~Fock and L.~Chekhov \cite{Foc, CheFoc1, CheFoc2} or, in a slightly different form, by R.~Kashaev \cite{Kash}; see also \cite{BonLiu, Liu1, GuoLiu}. This quantization takes advantage of the fact that, if one restricts to matrices with real coefficients,  a large subset  of $\RR_{\SL(\R)}(S)$ with non-empty interior  has a natural identification with the Teichm\"uller (or Fricke-Klein) space $\T(S)$, consisting of isotopy classes of all complete hyperbolic metrics on $S$. Thurston \cite{Thu} introduced for the Teichm\"uller space $\T(S)$ a set of coordinates, called shear coordinates, in which the Goldman-Weil-Petersson form is expressed in a particularly simple way. The \emph{quantum Teichm\"uller space} is a quantization of $\T(S)$ that is based on these shear coordinates. This construction requires the surface to have at least one puncture. 

A natural conjecture is that these two quantizations are ``essentially equivalent''. 

In the classical cases where $q=1$ and $A=\pm1$, the correspondence is relatively clear because of the identifications of $\mathcal S^{\pm1}(S)$ and $\widehat{\T}^1_S$ with algebras of functions on $\RR_{\SL(\C)}(S)$ and $\T(S)$. The only minor problem is that the functions considered in each case are not quite the same. 

The correspondence between the skein algebra $\mathcal S^{\pm1}(S)$ and the algebra of regular functions on $\RR_{\SL(\C)}(S)$ uses the trace functions $\Tr_K\col \RR_{\SL(\C)}(S) \to \R$, associated to all closed curves $K$ immersed in $S$, which to a homomorphism $r\col \pi_1(S) \to \SL(\C)$ associates the trace of $r(K) \in \SL(\C)$ (see \S\S \ref{sect:ClassicalTrace} and \ref{subsect:ClassicalSkeins} for technical details).  

Shear coordinates depend on the choice of some topological information, namely on the choice of an ideal triangulation $\lambda$ for the surface $S$. For a real representation $r_m\in\RR_{\SL(\R)}(S)$ corresponding to a hyperbolic metric $m\in  \T(S)$, the trace of $r_m(K)$ can then be explicitly computed (see \S\ref{subsect:ClassicalSkeins} for sign issues). This trace is actually expressed as a Laurent polynomial in the \emph{square roots} of the shear coordinates of $m$. This leads us to consider an algebra $\widehat{\mathcal Z}^1_\lambda$ consisting of rational fractions in the square roots of the shear coordinates, and to consider the algebra homomorphism
$$
\Tr_\lambda^1  \col \mathcal S^1(S) \to \widehat{\mathcal Z}^1_\lambda
$$
which to $[K]\in \mathcal S^1(S)$ associates the Laurent polynomial expressing the trace of $r_m(K)$ in terms of the shear coordinates of $m\in  \T(S)$. 

In the quantum case, one similarly introduces a non-commutative algebra $\ZZZ_\lambda$ consisting of rational fractions in certain skew-commuting variables associated to the square roots of the shear coordinates. When $q=\omega^4$, the quantum Teichm\"uller space $\widehat\T^q_S$ consists of those rational fractions in $\ZZZ_\lambda$ that involve only even powers of the variables. See \S \ref{sect:CFsquareRoot} for details.

\begin{thm}
\label{thm:MainThmIntro}
For $A=\omega^{-2}$, there  is an algebra homomorphism 
$$\Tr_\lambda^\omega  \col \SSS(S) \to \ZZZ_\lambda,$$
 depending continuously  on $\omega$ in an appropriate sense, which corresponds to the above homomorphism $\Tr_\lambda^1  \col \mathcal S^1(S) \to \mathcal Z^1_\lambda$ when $\omega=1$. In addition, the image $\Tr_\lambda^\omega([K]) \in \ZZZ_\lambda$ of every  $[K]\in \SSS(S)$ is a Laurent polynomial in the variables generating $\ZZZ_\lambda$. 
\end{thm}

The homomorphism $\Tr_\lambda^\omega$ is shown to be injective in Proposition~\ref{prop:QTraceInjective}. 

A major step in the construction of the quantum Teichm\"uller space $\widehat\T_S^q$ is to make it independent of a choice of ideal triangulation. The homomorphism $\Tr_\lambda^\omega$ of 
Theorem~\ref{thm:MainThmIntro} is similarly independent of choices. Making sense of this statement uses  work of C.~Hiatt in \cite{Hiatt} that extends to the square root set-up the original coordinate changes of Chekhov and Fock. More precisely, for any two ideal triangulations $\lambda$ and $\lambda'$ of the surface $S$, Hiatt constructs a coordinate change isomorphism $\Theta_{\lambda\lambda'}^\omega \col \ZZZ_{\lambda'}\to \ZZZ_\lambda$ that restricts to the identity on $\widehat\T^q_S$, considering the quantum Teichm\"uller space $\widehat\T^q_S$ as a subalgebra of both $ \ZZZ_{\lambda}$ and  $ \ZZZ_{\lambda'}$.

\begin{thm}
\label{thm:CoordinateChartInvarianceIntro}
Given two ideal triangulations $\lambda$ and $\lambda'$ of the surface $S$ and given an element $[K] \in \SSS(S)$ of the skein algebra of $S$, the coordinate change map
$$
\Theta_{\lambda\lambda'}^\omega \col \ZZZ_{\lambda'}\to \ZZZ_\lambda
$$
sends the Laurent polynomial $\Tr_{\lambda'}^\omega (K)$ to  the Laurent polynomial $\Tr_{\lambda}^\omega(K)$.
\end{thm}

While the proof of Theorem~\ref{thm:MainThmIntro} is rather elaborate, the proof of Theorem~\ref{thm:CoordinateChartInvarianceIntro} results from an easy application of the technology developed by Hiatt in \cite{Hiatt}.

Theorems~\ref{thm:MainThmIntro} and \ref{thm:CoordinateChartInvarianceIntro} were conjectured in \cite{Foc, CheFoc2}, and proved for certain small surfaces in \cite{ChePen1, Hiatt}. Our proof is much more 3--dimensional than these earlier attempts. The technical challenge is to figure out a ``good'' way to order the non-commuting variables in each monomial of the Laurent polynomials considered; this is a classical problem in mathematical physics, where it is known as the search for a quantum ordering. Our solution is based on a careful control of the elevations of the strands of a link $K$ in $S\times [0,1]$, with respect to the $[0,1]$ factor. The exposition that we give here is very computational, and involves a few miraculous identities (see in particular the proof of Proposition~\ref{prop:QuantumTraceWellDefined}) that the reader may find somewhat frustrating. Recent conversations with C.~Kassel seem to provide a more conceptual explanation for these identities, based on the fundamental representation of the dual $\SL(q)$ of the quantum group $\mathrm U_q( \mathrm{sl}_2)$; in particular, it might be possible to place our construction within the framework of \cite{BFK2, BFK3}. 

\medskip
The motivation for this work finds its origins in the respective advantages and drawbacks of the two points of view on the character variety $\RR_{\SL(\C)}(S)$, and in their impact on the corresponding quantizations. The algebraic geometric approach of $\RR_{\SL(\C)}(S)$, based on trace functions, is very natural and its coordinate functions use only polynomials; however, it is hard to extract much information from this description. Conversely, the shear coordinates for the Teichm\"uller space are very concrete and geometric, but they also are less intrinsic (in particular for hyperbolic surfaces with infinite area, for which additional data is needed), do not behave well under the operation of restriction to subsurfaces, and are not defined for closed surfaces. 
The same features can be found  at the quantum level. The skein algebra is very natural and occurs in many different contexts. However, its algebraic structure is quite difficult to handle at this point, except for small surfaces (see for instance \cite[\S3]{BonWon1} for a discussion). Conversely, the quantum Teichm\"uller space has a very simple algebraic structure (it is a quantum torus), but it suffers from the lack of canonicity inherited from the classical shear coordinates. 

One great advantage of the quantum Teichm\"uller space is that it has a very nice finite-dimensional representation theory, where an irreducible representation is essentially determined by a point in the  character variety $\RR_{\PSL(\C)}(S)$ \cite{BonLiu, BBL}. By composition with the trace homomorphism $\Tr_\lambda^\omega  \col \SSS(S) \to \ZZZ_\lambda$ provided by Theorem~\ref{thm:MainThmIntro}, one obtains a wide family of finite-dimensional representations of the skein algebra $ \SSS(S)$. These representations behave well with respect to the action of the mapping class group, and a great feature of the corresponding machinery is that it works even for closed surfaces \cite{BonWon1, BonWon2, BonWon3}. In particular, the results of the current paper represent a key technical step in a long-term program to study the representation theory of the skein algebra $ \SSS(S)$; see \cite{BonWon1} for a discussion.  

\medskip
\noindent\textbf{Acknowledgements.} We are grateful to Adam Sikora for helping us sorting out our ideas in the classical case where $A = \pm1$, and to Qingtao Chen for pointing out many misprints in the earlier versions of this paper. 

\section{The classical case}

\subsection{Ideal triangulations} The introduction was restricted to surfaces with no boundary, but it is convenient to allow boundary as well.

Let $S$ be an oriented  punctured surface with
boundary, obtained by removing finitely many points $v_1$,
$v_2$,
\dots, $v_p$ from a compact connected oriented surface
$\bar S$  with (possibly empty) boundary $\partial \bar S$. We require that each
component of $\partial \bar S$ contains at least one puncture $v_i$, that
there is at least one puncture, and that $\chi(S)<\frac d2$, where $\chi(S)$  is the Euler characteristic of $S$ and $d$ is the number of components of $\partial S$. These topological restrictions are
equivalent   to the existence of an \emph{ideal triangulation} for $S$,
namely a triangulation of the closed surface
$\bar S$ whose vertex set is exactly
$\{ v_1,
\dots, v_p\}$. In particular, an ideal triangulation
$\lambda$  has $n=-3\chi(S)+2d$ edges and $m=-2\chi(S)+d$ faces. Its edges
provide $-3\chi(S)+2d$ arcs $\lambda_1$, \dots,~$\lambda_n$ in $S$, going from
puncture to puncture, which decompose the surface $S$ into  $-2\chi(S)+d$
infinite triangles $T_1$, $T_2$, \dots, $T_m$ whose vertices sit ``at infinity'' at the punctures.
Note that $d$ of these $\lambda_i$ are just the boundary components
of $S$.

\subsection{The shear parameters}
Suppose that we are given a positive weight $X_i\in \R_+$ for each interior edge $\lambda_i$ of the ideal triangulation $\lambda$. We can associate to this data a group homomorphism $r\col \pi_1(S) \to \PSL(\R)$ as follows. 

Lift the ideal triangulation $\lambda$ to an ideal triangulation $\widetilde \lambda$ of the universal cover $\widetilde S$. We can then construct an orientation-preserving immersion $\widetilde f\col \widetilde S \to \HH^2$ from $\widetilde S$ to the hyperbolic plane $\HH^2$ such that:
\begin{enumerate}

\item $\widetilde f$ sends each face $\widetilde T$ of $\widetilde \lambda$ to an ideal triangle of $\HH^2$, delimited by three disjoint geodesics and touching the circle at infinity $\partial_\infty \HH^2$ in $3$ points;

\item when two faces $\widetilde T$ and $\widetilde T'$ meet along an edge $\widetilde \lambda_i$ that projects to the edge $\lambda_i$ of $\lambda$, then $\widetilde f( \widetilde T' )$ is obtained from $\widetilde f( \widetilde T)$ by performing a hyperbolic reflection across the geodesic $\widetilde f( \widetilde\lambda_i)$ followed by a hyperbolic translation of $\log X_i$ along the same geodesic $\widetilde f( \widetilde\lambda_i)$, if we orient $\widetilde f( \widetilde\lambda_i)$ by the boundary orientation of $\widetilde T$. 

\end{enumerate}

The immersion $\widetilde f$ is easily constructed stepwise, and uniquely determined up to isotopy of $\widetilde S$ respecting $\widetilde \lambda$, once we have chosen the image of a single face of $\widetilde S$. In particular, the family of the ideal triangles $\widetilde f(\widetilde T)\subset \HH^2$ is unique up to  an orientation-preserving isometry of $\HH^2$, namely up to composition by an element of $\PSL(\R)$. 

From the construction, it is immediate that there is a unique group homomorphism $r\col \pi_1(S) \to \PSL(\R)$ such that $\widetilde f(\gamma\widetilde T) = r(\gamma)(\widetilde T)$  for every face $\widetilde T$ of $\widetilde \lambda$. Since the family of the ideal triangles $\widetilde f(\widetilde T)$ is unique up to composition by an element of $\PSL(\R)$, $r$ is unique up to conjugation by an element of $\PSL(\R)$. 

We say that $r\col \pi_1(S) \to \PSL(\R)$ is associated to the \emph{shear parameters} $X_i \in \R_+$. 

\subsection{The classical trace function}
\label{sect:ClassicalTrace}

For a group homomorphism $r\col \pi_1(S) \to \PSL(\R)$ and immersion $\widetilde f \col \widetilde S \to \HH^2$ as above, consider a closed curve $ K $ immersed in $S$. 

The fact that $ K $ is immersed provides a natural lift $\widehat r(K)\in\SL(\R)$ of $r( K ) \in \PSL(\R)$. Indeed, lift $ K $ to an immersed path $\widetilde K  \col [0,1] \to \widetilde S$. Then $r( K ) \in \PSL(\R)$ is the unique orientation-preserving isometry of $\HH^2$ sending the point $\widetilde f \circ \widetilde  K  (0)\in \HH^2$ to $\widetilde f \circ \widetilde  K  (1)$, and sending the vector $(\widetilde f \circ \widetilde  K )' (0)$ to $(\widetilde f \circ \widetilde  K )' (1)$. Now, for every $t \in [0,1]$, we can consider the isometry $r( K )_t \in \PSL(\R)$ that sends the point $\widetilde f \circ \widetilde  K  (0)\in \HH^2$ to $\widetilde f \circ \widetilde  K  (t)$, and the vector $(\widetilde f \circ \widetilde  K )' (0)$ to a positive real multiple of $(\widetilde f \circ \widetilde  K )' (t)$. We now have a constructed a path $t\mapsto r( K )_t \in \PSL(\R)$ that joins $r( K )_0 = \Id_{\HH^2}$ to $r( K )_1 = r( K )$. This path defines an element of the universal cover of $\PSL(\R)$, which projects to an element $\widehat r(K)$ of the 2--fold cover  $\SL(\R)$ of $\PSL(\R)$. 

We are particularly interested in the trace $\Tr\, \widehat r(K)$ of $\widehat r(K) \in \SL(\R)$. Note that, when $K$ is just a small circle bounding a disk embedded in $S$, our designated lift $\widehat r(K)$  is minus the identity matrix of $\SL(\R)$, and $\Tr\, \widehat r(K)=-2$.

If the homomorphism $r\col \pi_1(S) \to \PSL(\R)$  is associated to shear parameters $X_i \in \R_+$ assigned to the edges of the ideal triangulation $\lambda$, the construction of the map $\widetilde f \col \widetilde S \to \HH^2$ and of the homomorphism $r$ is sufficiently explicit that $\widehat r(K) \in \SL(\R)$ can be explicitly computed. 

More precisely, suppose that $K$ transversely meets the edges $\lambda_{i_1}$, $\lambda_{i_2}$, \dots, $\lambda_{i_k}$, $\lambda_{i_{k+1}}=\lambda_{i_1}$, in this order. 
After crossing the edge $\lambda_{i_j}$, the curve $K$ enters a face $T$ of $\lambda$, which it exits through the edge $\lambda_{i_{j+1}}$. There are three possible choices for $\lambda_{i_{j+1}}$: it can be the edge immediately to the left as one enters $T$ through $\lambda_{i_j}$, the one immediately to the right, or it can be $\lambda_{i_j}$ again if $\gamma$ makes a U-turn in $T$. In addition, because $K$ is immersed, we can measure the amount by which the tangent to  $K$ turns between $\lambda_{i_j}$ and $\lambda_{i_{j+1}}$. We then define a matrix $M_j$ according to the various possible configurations.

If $\lambda_{i_{j+1}}$ is the edge immediately to the left as one enters $T$ through $\lambda_{i_j}$, let $t_j\in \Z$ denote the number of full turns to the left that  the tangent to $K$ makes between $\lambda_{i_j}$ and $\lambda_{i_{j+1}}$, and let $\epsilon_j = (-1)^{t_j} = \pm1$. Here the topological number of turns $t_j\in \Z$ is measured  so that $t_j=0$ when $K$ has no self-intersection between $\lambda_{i_j}$ and $\lambda_{i_{j+1}}$; in fact, $t_j$ has the same parity as the number of double points of $K$ between $\lambda_{i_j}$ and $\lambda_{i_{j+1}}$. In this case, define
$$
M_j =  \begin{pmatrix}
\epsilon_j & \epsilon_j \\0& \epsilon_j
\end{pmatrix}.
$$

For the analogous case where $\lambda_{i_{j+1}}$ is the edge immediately to the right as one enters $T$ through $\lambda_{i_j}$, let again $t_j\in \Z$ denote the number of full turns to the left that  the tangent to  $K$ makes between $\lambda_{i_j}$ and $\lambda_{i_{j+1}}$, and set $\epsilon_j = (-1)^{t_j} = \pm1$. Then define
$$
M_j =  \begin{pmatrix}
\epsilon_j & 0 \\ \epsilon_j  & \epsilon_j
\end{pmatrix}.
$$

In the case of a U-turn, where $\lambda_{i_{j+1}} = \lambda_{i_j}$, let $t_j\in \Z$ be defined so that  the tangent to  $K$ makes $2t_j +1$ half-turns to the left between $\lambda_{i_j}$ and $\lambda_{i_{j+1}}$, and set again $\epsilon_j = (-1)^{t_j} = \pm1$. Then define
$$
M_j =  \begin{pmatrix}
0& \epsilon_j  \\ -\epsilon_j  & 0
\end{pmatrix}.
$$

Finally, having defined $M_j$ in every case,  consider for $X >0$ the matrix
$$
S(X)=
\begin{pmatrix}
X^{\frac12}&0\\0& X^{-\frac12}
\end{pmatrix}.
$$

\begin{lem} 
\label{lem:PlusMinusTrace}
Up to conjugation by an element of $\SL(\R)$,
$$
\widehat r(K) = S(X_{i_1}) M_1 S(X_{i_2}) M_2 \dots S(X_{i_k}) M_k
$$
where the matrices $M_j$ and $S(X_{i_j})$ are associated as above to the way the immersed curve $K$ crosses the edges of the ideal triangulation $\lambda$, and where $X_i \in \R_+$ are the shear parameters defining the homomorphism $r\col \pi_1(S) \to \PSL(\R)$. 
\end{lem}
\begin{proof}
This is an easy exercise in hyperbolic geometry. See for instance Exercises~8.5--8.7 and 10.14 in \cite{BonBook}.
\end{proof}

\subsection{State sums}

As preparation for the quantum extension, we now give a state sum formula for the trace $\Tr\, \widehat r(K)$ of the above element $\widehat r(K) \in \SL(\R)$. 

Let a \emph{state} assign a sign $s_1$, $s_2$, \dots, $s_k$, $s_{k+1}=s_1\in \{ +, -\}$ to each point where $K$ crosses an edge $\lambda_{i_j}$ of $\lambda$, in this order. For $j=1$, $2$, \dots, $k$, write the matrix $M_j$ defined above as 
$$
M_j = 
\begin{pmatrix}
m_j^{++} & m_j^{+-} \\
m_j^{-+}  & m_j^{--}
\end{pmatrix}
$$
with $m_j^{\pm\pm} =0$, $+1$ or $-1$.  
Then, the following formula immediately follows from Lemma~\ref{lem:PlusMinusTrace} combined with elementary linear algebra.

\begin{lem}
\label{lem:ClassicalStateSum}
$$
\Tr \, \widehat r(K) = \sum_s 
  m_1^{s_1 s_2}
 m_2^{s_2 s_3}
 \dots
  m_k^{s_k s_1}
  X_{i_1}^{\frac {s_1}2} 
 X_{i_2}^{\frac {s_2} 2} 
 \dots
  X_{i_k}^{\frac {s_k}2} 
$$
where the sum is over all possible states $s$ for $K$ and $\lambda$, and where in the exponents we identity the sign $s_j = \pm$ to the number $s_j=\pm1$. \qed
\end{lem}

\section{The quantum Teichm\"uller space}
\label{sect:QuantumTeich}

\subsection{The Chekhov-Fock algebra of an ideal triangulation}

Let $T_1$, $T_2$, \dots, $T_m$ be the faces of the ideal triangulation $\lambda$. Index the sides of each face $T_j$ as $\lambda_{j1}$, $\lambda_{j2}$, $\lambda_{j3}$, in such a way that they occur in this order clockwise around $T_j$. We then  associate to $T_j$ a copy $\T_{T_j}^q$ of the \emph{triangle algebra}, generated by three elements $X_{j1}$, $X_{j2}$, $X_{j3}$ and their inverses $X_{j1}^{-1}$, $X_{j2}^{-1}$, $X_{j3}^{-1}$, and defined by the relations that $X_{j1}X_{j2}=q^2 X_{j2}X_{j1}$, $X_{j2}X_{j3}=q^2 X_{j3}X_{j2}$ and $X_{j3}X_{j1}=q^2 X_{j1}X_{j3}$. We here think of each generator $X_{j a}$ as being associated to the side $\lambda_{ja}$ of $T_j$. 

In the tensor product algebra $\T_{T_1}^q \otimes  \dots \otimes \T_{T_m}^q = \bigotimes_{j=1}^m  \T_{T_j}^q$, we now associate to the edge $\lambda_i$ of $\lambda$ an element $X_i$, defined by:
\begin{enumerate}
\item $X_i = X_{ja} \otimes X_{kb}$ if $\lambda_i$ separates two
distinct faces
$T_j$ and
$T_k$, and if
$X_{ja}\in \T_{T_{j}}^q$ and $X_{kb}\in \mathcal
T_{T_{k}}^q$ are the generators associated to the sides of
$T_j$ and $T_k$ corresponding to
$\lambda_i$;
\item $X_i = q^{-1} X_{ja} X_{jb} =q X_{jb} X_{ja}$ if $\lambda_i$
corresponds to two sides of the same face $T_j$, if $X_{ja}$,
$X_{jb}\in\T_{T_j}^q$ are the generators associated to these two
sides, and if $X_{ja}$ is associated to the side that comes first when going
counterclockwise around their common vertex. 
\end{enumerate}
By convention, when describing an element $Z_1 \otimes  \dots \otimes Z_m$ of $\bigotimes_{j=1}^m  \T_{T_j}^q$, we omit in the tensor product those $Z_j$ that are equal to the identity element $1$ of $\T_{T_j}^q$.

The \emph{Chekhov-Fock algebra} of the ideal triangulation $\lambda$ is the subalgebra $\T_\lambda^q$ of $\bigotimes_{j=1}^m  \T_{T_j}^q$ generated by the elements $X_i$ associated as above to the edges of $\lambda$, and by their inverses $X_i^{-1}$. 

Note that $X_iX_j = q^{2\sigma_{ij}} X_jX_i$  where the integers $\sigma_{ij}\in \{ 0, \pm1, \pm2\}$ are
defined as follows: Let $a_{ij}$ be the number of angular sectors  delimited  by $\lambda_i$ and  $\lambda_j$ in the
faces of
$\lambda$, and with $\lambda_i$
coming first counterclockwise; then $\sigma_{ij} = a_{ij}-a_{ji}$. 

\subsection{Coordinate changes between Chekhov-Fock algebras}

As one switches from one ideal triangulation $\lambda$ to another ideal triangulation $\lambda'$, the geometry of the Teichm\"uller space provides coordinate changes between the shear coordinates associated to $\lambda$ and those associated to $\lambda'$. Because shear coordinates can be expressed as cross-ratios, one easily sees that these coordinate changes are given by rational maps. 

In the quantum case, there is no underlying geometry to provide us with similar coordinate changes, and one has to find algebraic isomorphisms that have the required properties. 

As in the classical case, these will involve rational fractions, and we consequently have to introduce the
\emph{fraction
division algebra} 
$\widehat{\T}_\lambda^q$ of the Chekhov-Fock algebra $\T_\lambda^q$. Such a fraction division algebra exists because $\T_\lambda^q$ satisfies the so-called Ore Condition; see for instance
\cite{Coh, Kass}. In practice,
$\widehat{\T}_\lambda^q = \C (X_1,
\dots, X_n)_\lambda^q$ consists of
non-commutative rational fractions in the
variables $X_1$, \dots, $X_n$ which are manipulated
according to the $q$-commutativity relations
$X_iX_j = q^{2\sigma_{ij}} X_jX_i$. 

L. Chekhov and V. Fock \cite{Foc, CheFoc1, CheFoc2} (and R. Kashaev \cite{Kash} in the context of length coordinates) construct such coordinate isomorphisms; see also \cite{BonLiu, Liu1}.

\begin{thm}[Chekhov-Fock, Kashaev]
\label{thm:CheFockCoordinateChanges}
There exists a family of algebra isomorphisms
$$ \Phi_{\lambda\lambda'}^q \col \widehat{\T}_{\lambda'}^q
 \to \widehat{\T}_\lambda^q,
$$
defined for any two ideal triangulations $\lambda$, $\lambda'$, such that
$$\Phi_{\lambda\lambda''}^q =\Phi_{\lambda\lambda'}^q \circ \Phi_{\lambda'\lambda''}^q $$ 
for any three ideal triangulations $\lambda$, $\lambda'$ and $\lambda''$. \qed
\end{thm}

This enables us to define the  \emph{quantum Teichm\"uller space} $\widehat{\T}_S^q$ of the punctured surface $S$  as the quotient
$$
\widehat{\T}_S^q = \coprod_\lambda \widehat{\T}_\lambda^q /\sim 
$$
of the disjoint union of the $ \widehat{\T}_\lambda^q$ of all ideal triangulations $\lambda$ of $S$, where the equivalence relation $\sim$ identifies $\widehat{\T}_\lambda^q $  to 
$\widehat{\T}_{\lambda'}^q$ by the coordinate change isomorphism $\Phi_{\lambda\lambda'}^q$. Note that the property that $\Phi_{\lambda\lambda''}^q =\Phi_{\lambda\lambda'}^q \circ \Phi_{\lambda'\lambda''}^q $  is crucial to guarantee that $\sim$ is an equivalence relation. This property is much stronger than one could have thought at first glance, as indicated by the uniqueness result of \cite{Bai}. 

Because the $\Phi_{\lambda\lambda'}^q$ are algebra isomorphisms, the quantum Teichm\"uller space  $\widehat{\T}_S^q$ inherits an algebra structure from the $\widehat{\T}_\lambda^q $.

\subsection{The Chekhov-Fock square root algebra}
\label{sect:CFsquareRoot}

The formulas of Lemma~\ref{lem:ClassicalStateSum} involve \emph{square roots} of shear coordinates. This lead us to consider formal square roots $Z_i = X_i^{\frac12}$ of the generators of the Chekhov-Fock algebra $\T_\lambda^q$. 

In practice, one just considers the Chekhov-Fock algebra $\T_\lambda^\omega$ associated to a fourth root $\omega= q^{\frac14}$ of $q$. To avoid confusion, we denote by $Z_i$ the generator of $\T_\lambda^\omega$ associated to the edge $\lambda_i$ of $\lambda$ while, as before, $X_i$ is the generator of $\T_\lambda^q$ associated to the same $\lambda_i$. Then,  there exists an injective algebra homomorphism $\T_\lambda^q \to \T_\lambda^\omega$  associating the element $Z_i^2$ to the generator $X_i$, so that we can consider $\T_\lambda^q $ as a subalgebra of $ \T_\lambda^\omega$. This also induces a similar inclusion $\widehat{\T}_\lambda^q \subset \widehat{\T}_\lambda^\omega$ between the corresponding fraction division algebras. 

In the classical case, the coordinate changes between square roots of shear coordinates are not as nice as those between shear coordinates, because they are not rational anymore. The same consequently holds in the quantum setup. However, there is a subalgebra of the  algebra $\T_\lambda^\omega$ which is better behaved with respect to coordinate changes.

A monomial $Z_1^{k_1}Z_2^{k_2}\dots Z_n^{k_n}$ in the generators $Z_i$ of $\T_\lambda^\omega$ is said to be \emph{balanced} if, for every triangle face $T_j$ of $\lambda$,  the exponents $k_i$ of the generators $Z_i$ associated to the three sides of $T_j$ add up to an even number. (When the same edge $\lambda_i$ corresponds to two distinct sides of $T_j$, the exponent $k_i$ is counted twice in the sum.) This is equivalent to the property that there exists a homology class $\alpha \in H_1(S; \Z_2)$ such that the class of the exponent $k_i$ in $\Z_2$ is equal to the algebraic intersection number of $\alpha$ with the edge $\lambda_i$. In this case, we will say that the monomial $Z_1^{k_1}Z_2^{k_2}\dots Z_n^{k_n}$ is $\alpha$--\emph{balanced}. 

In the Chekhov-Fock algebra $\T_\lambda^\omega$, let $\ZZ_\lambda$ denote the linear subspace generated by all balanced monomials. Note that it splits as a direct sum
$$
\ZZ_\lambda = \bigoplus_{\alpha\in H_1(S, \Z_2)} \ZZ_\lambda(\alpha)
$$
where $\ZZ_\lambda(\alpha)$ denotes the linear subspace generated by all $\alpha$--balanced monomials. The product of an element of $\ZZ_\lambda(\alpha)$ with an element of $\ZZ_\lambda(\beta )$ belongs to $\ZZ_\lambda(\alpha+\beta)$, so that $\ZZ_\lambda$ is a subalgebra of $\T_\lambda^\omega$. We will refer to $\ZZ_\lambda$ as the \emph{Chekhov-Fock square root algebra} of the ideal triangulation $\lambda$. 

Note that $\ZZ_\lambda(0)$ coincides with the subalgebra $\T_\lambda^q \subset \T_\lambda^\omega$ generated by the $Z_i^2 = X_i$. 

In the fraction algebra $\widehat{\T}_\lambda^\omega$, we now consider the subset $\ZZZ_\lambda$ consisting of all fractions $P/Q$ where $P\in \ZZ_\lambda$ and $Q\in \ZZ_\lambda(0)= \T_{\lambda}^q$. One easily checks that $\ZZZ_\lambda$ is a subalgebra of $\widehat{\T}_\lambda^\omega$. It also contains the Chekhov-Fock fraction algebra $ \widehat{\T}_\lambda^q$.

 In \cite{Hiatt}, Chris Hiatt extends  the Chekhov-Fock coordinate change isomorphism $\Phi_{\lambda\lambda'}\col  \widehat{\T}_{\lambda'}^q \to  \widehat{\T}_{\lambda}^q$ of Theorem~\ref{thm:CheFockCoordinateChanges} to the Chekhov-Fock square root algebra $\ZZ_\lambda$.

\begin{thm}[Hiatt]
\label{thm:SquareRootCoordChanges}
When $q=\omega^4$, there exists for any two ideal triangulations $\lambda$, $\lambda'$ an algebra isomorphism
$$
\Theta_{\lambda\lambda'}^\omega \col \ZZZ_{\lambda'} \to \ZZZ_\lambda
$$
extending the Chekhov-Fock coordinate change isomorphism $\Phi_{\lambda\lambda'}\col  \widehat{\T}_{\lambda'}^q \to  \widehat{\T}_{\lambda}^q$. In addition, $\Theta_{\lambda\lambda''}^\omega= \Theta_{\lambda\lambda'}^\omega\circ  \Theta_{\lambda'\lambda''}^\omega$ for any three ideal triangulations $\lambda$, $\lambda'$ and $\lambda''$. 
\end{thm}

\begin{proof}

Hiatt does not quite prove the result in this form, so we need to explain how to obtain it from \cite[\S\S 6--7]{Hiatt}. 

For every $\alpha\in H_1(S; \Z_2)$, let $\ZZZ_\lambda(\alpha)$ consist of all fractions $P/Q$ where $P\in \ZZ_\lambda(\alpha)$ and $Q\in \ZZ_\lambda(0)= \T_{\lambda}^q$, so that 
$$
\ZZZ_\lambda = \bigoplus_{\alpha\in H_1(S, \Z_2)} \ZZZ_\lambda(\alpha).
$$

Hiatt constructs in \cite[\S6]{Hiatt} a linear map $
\Theta_{\lambda\lambda'}^\omega \col \ZZZ_{\lambda'}(\alpha) \to \ZZZ_\lambda(\alpha)
$. The construction of this map  in \cite{Hiatt} \emph{a priori} depends on the choice of an 1--dimensional submanifold $K$ immersed in $ S$ and representing the homology class $\alpha \in H_1(S; \Z_2)$. However, it easily follows from \cite[Lemma~17]{Hiatt} that this map depends  only on $\alpha$. 

Linearly extend these maps $
\Theta_{\lambda\lambda'}^\omega \col \ZZZ_{\lambda'}(\alpha) \to \ZZZ_\lambda(\alpha)
$ to a linear map $
\Theta_{\lambda\lambda'}^\omega \col \ZZZ_{\lambda'} \to \ZZZ_\lambda
$. 

To show that this is an algebra homomorphism we need to check that, for every $A\in \ZZZ_{\lambda'}(\alpha)$ and $B\in \ZZZ_{\lambda'}(\beta)$, $\Theta_{\lambda\lambda'}^\omega(A)\Theta_{\lambda\lambda'}^\omega(B)= \Theta_{\lambda\lambda'}^\omega(AB)$ in $\ZZZ_\lambda(\alpha+\beta)$. This is an immediate consequence of \cite[Sublemma~19]{Hiatt} and of the construction of $\Theta_{\lambda\lambda'}^\omega$. 

The fact that the restriction of $\Theta_{\lambda\lambda'}^\omega$ to $\ZZZ_{\lambda'}(0) = \widehat{\T}_{\lambda'}^q$ coincides with $\Phi_{\lambda\lambda'}^q$ immediately follows from  its construction in \cite{Hiatt}. The property that $\Theta_{\lambda\lambda''}^\omega= \Theta_{\lambda\lambda'}^\omega\circ  \Theta_{\lambda'\lambda''}^\omega$ is proved in \cite[Theorem~25]{Hiatt}. 
\end{proof}

\begin{rem}
For $A \in \ZZZ_{\lambda}$, the operator point of view of \cite{Foc, CheFoc1, CheFoc2} much more easily provides a natural square root $\Theta_{\lambda\lambda'}^\omega(A)$ of $\Phi_{\lambda\lambda'}^q(A^2)$. The real content of Theorem~\ref{thm:SquareRootCoordChanges} is that this square root can be expressed as a rational fraction in the generators $Z_i=X_i^{\frac12}$. The restriction to $\ZZZ_{\lambda'}$ is here crucial. 
\end{rem}

\section{The skein algebra}

\subsection{Links and skeins}\label{sect:skein}

We begin with the \emph{framed link algebra} $\mathcal K(S)$. This is the vector space (over $\C$, say) freely generated by the isotopy classes of (unoriented) 1--dimensional framed submanifolds $K \subset S \times [0,1]$ such that:
\begin{enumerate}
\item $\partial K =  K \cap \partial(S \times [0,1])$ consists of finitely many points in $(\partial S) \times [0,1]$;
\item at every point of $\partial K$, the framing is vertical, namely parallel to the $[0,1]$ factor, and  points in the direction of $1$;
\item for every component $k$ of $\partial S$, the points of $\partial K$ that are in $k\times [0,1]$ sit at different elevations, namely have different $[0,1]$--coordinates. 
\end{enumerate}
An isotopy of such framed submanifolds of course is required to respect all three conditions. 
The third condition will turn out to be crucial for our analysis.

Perhaps we should have begun by specifying what we mean by a \emph{framing} for $K$. For us here, a framing  is a continuous choice of a vector transverse to $K$ at each point of $K$.  

The vector space $\mathcal K(S)$ can be endowed with a multiplication, where the product of $K_1$ and $K_2$ is defined by the framed link $K \subset S\times[0,1]$ that is the union of $K_1$ rescaled in $S\times [0, \frac12]$ and $K_2$ rescaled in $S\times [\frac12, 1]$. In other words, the product $K_1K_2$ is defined by superposition of the framed links $K_1$ and $K_2$. 
Note that this \emph{superposition operation} is compatible with isotopies, and therefore provides a well-defined algebra structure on $\mathcal K(S)$. 

\begin{figure}[htbp]

\SetLabels
( .5 * -.4 ) $K_0$ \\
( .1 * -.4 )  $K_1$\\
(  .9*  -.4) $K_\infty$ \\
\endSetLabels
\centerline{\AffixLabels{\includegraphics{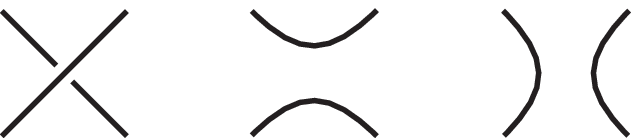}}}
\vskip 10pt
$$
 [K_1] = A^{-1}[K_0] + A [K_\infty]
 $$
\caption{The skein relation}
\label{fig:skein}
\end{figure}

Three links $K_1$, $K_0$ and $K_\infty$ in $S\times[0,1]$ form a \emph{Kauffman triple} if the only place where they differ is above a small disk in $S$, where they are as represented in Figure~\ref{fig:skein} (as seen from above) and where the framing is vertical and pointing upwards (namely the framing is parallel to the $[0,1]$ factor and points towards $1$).

The \emph{Kauffman skein algebra} $\SSS(S)$ is the quotient of the framed link algebra $\mathcal K(K)$ by the two-sided ideal generated by all elements $K_1 - A^{-1}K_0 - A K_\infty$ as $(K_1, K_0, K_\infty)$ ranges over all Kauffman triples. The superposition operation descends to a multiplication in $\SSS(S)$, endowing $\SSS(S)$ with the structure of an algebra. The class $[\varnothing]$ of the empty link is an identity element in $\SSS(S)$, and is usually denoted by $1$.

 An element $[K]\in \SSS(S)$, represented by a framed link $K \in \mathcal K(S)$, is a \emph{skein} in $S$. The construction is defined to ensure that the \emph{skein relation}
 $$
 [K_1] = A^{-1}[K_0] + A [K_\infty]
 $$
 holds in $\SSS(S)$ for every Kauffman triple $(K_1, K_0, K_\infty)$. 

\subsection{The classical cases $\mathbf{A=\pm 1}$}
\label{subsect:ClassicalSkeins}

In \cite{Bull1, Bull2, BFK1, PrzS}, D.~Bullock, C.~Frohman, J.~Kania-Bartoszy\'nska, J.~Przytycki and A.~Sikora observe that there is a strong connection between the skein algebra with $A=-1$ and the character  variety
$$
\RR_{\SL(\C)}(S) = \{ \text{group homomorphisms } r\col \pi_1(S) \to \SL(\C) \} \db \SL(\C).
$$
The quotient is under the action by conjugation, and should be understood in the sense of geometric invariant theory \cite{Mum} to avoid pathologies near the reducible homomorphisms. 

Note that, for every $A$, there is a unique algebra homomorphism $\SSS(S) \to \C$ that sends each non-empty skein $[K]\in \SSS(S)$ to $0$ and sends the empty skein $[\varnothing]=1$ to $1$. This homomorphism is the \emph{trivial homomorphism} $\SSS(S) \to \C$.

\begin{thm}[\cite{Bull1, Bull2, BFK1,  PrzS}]
\label{thm:Bullock}
Assume that the surface $S$ has no boundary (but is still allowed to have punctures), and consider the skein algebra $\mathcal S^{-1}(S) $ corresponding to $A=-1$. Every group homomorphism  $ r\col \pi_1(S) \to \SL(\C)$ defines a unique non-trivial algebra homomorphism $T_r \col \mathcal S^{-1}(S) \to \C$ by the property that 
$$
T_r( [K] ) = -\Tr \, r(K)
$$
for every connected skein $[K]\in \mathcal S^{-1}(S)$. 

Conversely, every  non-trivial algebra homomorphism $T \col \mathcal S^{-1}(S) \to \C$ is associated to a unique $r\in \RR_{\SL(\C)}(S)$ in this way. \qed
\end{thm}

Note that the definition of $r(K)\in \SL(\C)$ implicitly supposes the choice of an orientation for the closed curve $K$. However, reversing this orientation replaces $r(K)$ by its inverse, and leaves the trace $\Tr \, r(K)$ unchanged.

There is a similar result for the other case where the skein algebra $\SSS(S)$ is commutative, corresponding to $A=1$. This statement uses the correspondence $\SSS(S) \cong \mathcal S^{-A}(S)$ established by J.~Barrett \cite{Barr}, and requires the use of spin structures. 

Let $\Spin(S)$ be the set of isotopy classes of spin structures on $S$ or, equivalently, the set of isotopy classes of spin structures on $S\times[0,1]$. Any two elements of $\Spin(S)$ differ by an obstruction in $H^1(S; \Z_2)$, so that there is an action of $H^1(S; \Z_2)$ on $\Spin(S)$. 

Similarly, the cohomology group $H^1(S; \Z_2)$ acts on $\RR_{\SL(\C)}(S)$ by the property that, if $\alpha \in H^1(S;\Z_2)$ and $r\col \pi_1(S) \to \SL(\C)$, then $\alpha r\in \RR_{\SL(\C)}(S)$ associates $(-1)^{\alpha(\gamma)} r(\gamma) \in \SL(\C)$ to $\gamma\in \pi_1(S)$. Note that the quotient of $\RR_{\SL(\C)}(S)$ under this action of $H^1(S; \Z_2)$ is just the character variety
$$
\RR_{\PSL(\C)} (S) = \{ \text{group homomorphisms } r\col \pi_1(S) \to \PSL(\C) \} \db \PSL(\C)
$$
of homomorphisms valued in $\PSL(\C)$ instead of $\SL(C)$. (It is here important that $S$ is non-compact so that, because the fundamental $\pi_1(S)$ is free, every homomorphism $\pi_1 \to \PSL(\C)$ lifts to $\SL(\C)$.)

We can then combine these actions of $H^1(S; \Z_2)$  on $\Spin(S)$ and $\RR_{\SL(\C)}(S)$, and consider the twisted product
$$
\RS(S)= \RR_{\PSL(\C)}(S) \widetilde\times\, \Spin(S) = \bigl(\RR_{\SL(\C)}(S) \times \Spin(S) \bigr)/H^1(S; \Z_2).
$$
Note that, just like $\RR_{\SL(\C)}(S)$, this twisted product $\RS(S)$ is a finite cover of $\RR_{\PSL(\C)}(S)$ with fiber $H^1(S;\Z_2) \cong \Spin(S)$. 

If $\sigma\in \Spin(S)$ is a spin structure and $K$ is a framed knot in $S\times[0,1]$, the monodromy of the framing of $K$ with respect to $\sigma$ defines an element $\sigma(K) \in \Z_2$. If, in addition, we are given a group homomorphism $r\col \pi_1(S) \to \SL(\C)$, we can consider the element $T_{(r, \sigma)}(K) = (-1)^{\sigma(K)} \Tr\, r(K)$. Note that $T_{(r, \sigma)}(K)$ is invariant under the action of $H^1(S;\Z_2)$ on the pair $(r, \sigma)$, and therefore depends only on the class of $(r,\sigma)$ in $\RS(S)$.

\begin{thm}
\label{thm:BullockPSL}
Assume that the surface $S$ has no boundary (but is still allowed to have punctures), and consider the skein algebra $\mathcal S^{1}(S) $ corresponding to $A=+1$. Every group homomorphism  $ r\col \pi_1(S) \to \SL(\C)$ and spin structure $\sigma\in \Spin(S)$ define a unique algebra homomorphism $T_{(r, \sigma)} \col \mathcal S^{1}(S) \to \C$ by the property that 
$$
T_{(r, \sigma)}( [K] ) = (-1)^{\sigma(K)} \Tr \, r(K)
$$
for every connected skein $[K]\in \mathcal S^{1}(S)$. This homomorphism $T_{(r, \sigma)}$ is  non-trivial  and depends only on the class of $(r,\sigma)$ in $\RS(S)$. 

Conversely, every  non-trivial  algebra homomorphism $T \col \mathcal S^{1}(S) \to \C$ is associated to a unique element $(r, \sigma) \in \RS(S) $  in this way.
\end{thm}
\begin{proof}
Fix a spin structure $\sigma\in \Spin(S)$. Then Barrett \cite{Barr} defines an algebra isomorphism $\mathcal S^1(S) \to \mathcal S^{-1}(S)$ associating $(-1)^{k+\sigma(K)}[K] \in \mathcal S^{-1}(S)$ to every skein $[K] \in \mathcal S^1(S)$ represented by a link $K$ with $k$ components (see \cite[\S 2]{PrzS} for a proof that this is an algebra homomorphism). The result then immediately follows by combining Theorem~\ref{thm:Bullock} with  this correspondence. 
\end{proof}

To connect the set-up of \S \ref{sect:ClassicalTrace} to Theorem~\ref{thm:BullockPSL}, consider a hyperbolic metric $m\in \T(S)$. It is convenient to move to a 3--dimensional framework, by extending  $m \in \T(S)$ to a 3--dimensional hyperbolic metric on a small thickening  $S \times (0,1)$ of $S = S \times \{ \frac12\}$. We can even consider the more general case of a hyperbolic metric $m$ on  $S \times (0,1)$, not necessarily complete. Classically,  this hyperbolic metric $m$  on  $S \times (0,1)$ has a well-defined monodromy homomorphism $r\in \RR_{\PSL(\C)}(S)$. What seems less well-known is that $m$ provides additional spin information, and uniquely determines an element of the twisted product $\RS(S)= \RR_{\PSL(\C)}(S) \widetilde\times\, \Spin(S) $.

A spin structure $\sigma\in \Spin(S)$ specifies a way to lift the monodromy homomorphism $r \col \pi_1(S) \to \PSL(\C)$ to a homomorphism $r _\sigma\col \pi_1(S) \to \SL(\C)$ as follows. For this, first extend $\sigma$ to a spin structure on the thickened surface $S\times (0,1)$. Then consider a  \emph{developing map}  for the metric $m$, namely an isometric immersion $\widetilde f \col \widetilde S \times (0,1) \to \HH^3$ from the universal cover $\widetilde S \times (0,1)$ to the hyperbolic space $\HH^3$ that is equivariant with respect to the monodromy $r \col \pi_1(S) \to \PSL(\C)$. 

Represent $[\gamma]\in \pi_1(S)= \pi_1(S; x_0)$ by a path $\gamma \col [0,1] \to S\times(0,1)$ with $\gamma(0)=\gamma(1) = x_0$. Pick an arbitrary orthonormal frame $F(t)$ at each $\gamma(t)$, depending continuously on $t$ and such that $F(1)=F(0)$ at  $\gamma(0)=\gamma(1) = x_0$. Lift $\gamma$ to $\widetilde\gamma \col [0,1] \to \widetilde S \times(0,1)$, and $F(t)$ to an orthonormal frame $\widetilde F(t)$ at $\widetilde\gamma(t)$. For every $t\in [0,1]$, we can now consider the unique isometry $r(\gamma)_t\in \PSL(\C)$ of $\HH^3$ that sends $\widetilde f(\widetilde\gamma(0))$ to $\widetilde f(\widetilde\gamma(t))$ and $\widetilde f(\widetilde F(0))$ to $\widetilde f (\widetilde F(t))$. By construction, $r(\gamma)_0 = \Id_{\HH^3}$ and $r(\gamma)_1=r(\gamma)$. The path $t\mapsto r(\gamma)_t$ then defines a lift $\widehat r(\gamma)\in \SL(\C)$ of $r(\gamma)=r(\gamma)_1$ to the universal cover $\SL(\C)$ of $\PSL(\C)$. This $\widehat r(\gamma)\in \SL(\C)$ clearly depends on the framing $F$ but, if $\sigma(F)\in \Z_2$ denotes the monodromy of the framing $F$ around $\gamma$ with respect to the spin structure $\sigma$, 
$$
r_\sigma (\gamma) = (-1)^{\sigma(F)} \widehat r(\gamma) \in \SL(\C)
$$
does not. One easily checks that this defines a group homomorphism $r_\sigma \col \pi_1(S) \to \SL(\C)$. 

By definition of the action of $H^1(S; \Z_2)$ on $ \RR_{\PSL(\C)}(S)$ and $\Spin(S)$, a different choice of spin structure $\sigma \in \Spin(S)$ does not change the class of $(r_\sigma, \sigma)$ in $\RS(S) = \bigl(\RR_{\SL(\C)}(S) \times \Spin(S) \bigr)/H^1(S; \Z_2)$. 

\begin{prop}
The element $(r_\sigma, \sigma)$ in $\RS(S) = \RR_{\PSL(\C)}(S) \widetilde\times\, \Spin(S)$ depends only on the hyperbolic metric $m$ on $S\times (0,1)$.  \qed
\end{prop}

For a connected skein $[K]\in \mathcal S^1(S)$, note that the element $T_{(r_\sigma, \sigma)}([K])$ associated to $(r_\sigma, \sigma) \in \RS(S)$ by Theorem~\ref{thm:BullockPSL} is just the trace of $\widehat r(K)$ defined above. As a consequence, in the situation of \S \ref{sect:ClassicalTrace} where $m$ comes from a 2-dimensional hyperbolic metric on $S$, where the link $K$ is obtained by perturbing a curve immersed in $S$ to remove double points, and where the framing is chosen everywhere vertical, $T_{(r_\sigma, \sigma)}([K])$ is exactly the trace $\Tr\, \widehat r(K)$ considered in \S \ref{sect:ClassicalTrace}.

\subsection{Gluing skeins}

In addition to the multiplication by superposition, there is another operation which can be performed on framed links and skeins. 

Given two surfaces $S_1$ and $S_2$ and two boundary components $k_1 \subset \partial S_1$ and $k_2 \subset \partial S_2$, we can glue $S_1$ and $S_2$ by identifying $k_1$ and $k_2$ to obtain a new oriented  surface $S$. There is a unique way to perform this gluing so that the orientations of $S_1$ and $S_2$ match to give an orientation of $S$. We allow the ``self-gluing'' case, where the surfaces $S_1$ and $S_2$ are equal as long as the boundary components $k_1$ and $k_2$ are distinct. If we are given an ideal triangulation $\lambda_1$ of $S_1$ and an ideal triangulation $\lambda_2$ of $S_2$, these two triangulations fit together to give an ideal triangulation $\lambda$ of the glued surface $S$.

Now, suppose in addition that we are given skeins $[K_1] \in \SSS(S_1)$ and $[K_2]\in\SSS(S_2)$ such that $K_1 \cap( k_1 \times[0,1])$ and $K_2 \cap( k_2 \times[0,1])$ have the same number of points. We can then arrange by an isotopy of framed links that $K_1$ and $K_2$ fit together to give a framed link $K\subset  S \times [0,1]$; note that it is here important that the framings be vertical pointing upwards on the boundary, so that they fit together to give a framing of $K$. By our hypothesis that the points of $K_1 \cap( k_1 \times[0,1])$ (and of $K_2 \cap( k_2 \times[0,1])$ sit at different elevations, the framed link $K$ is now uniquely determined up to isotopy. Also, this operation is well behaved with respect to the skein relations, so that $K$ represents a well-defined element $[K] \in \SSS(S)$. 
We will say that $[K] \in \SSS(S)$ is \emph{obtained by gluing} the two skeins $[K_1] \in \SSS(S_1)$ and $[K_2]\in\SSS(S_2)$.

\subsection{The main theorem}
\label{subsect:MainThm}

Let a \emph{state} for a skein $[K]\in \SSS(S)$ be the assignment $s\col \partial K \to \{ +, -\}$ of a sign $\pm$ to each point of $\partial K$. Let $\SSSS(S)$ be the algebra consisting of linear combinations of \emph{stated skeins}, namely of  skeins endowed with states. 

In the case where $K \in \SSS(S)$ is obtained by gluing the two skeins $K_1 \in \SSS(S_1)$ and $K_2\in\SSS(S_2)$, the states $s\col \partial K \to \{ +,-\}$, $s_1\col \partial K_1 \to \{ +,-\}$, $s_2\col \partial K_2 \to \{ +,-\}$ are \emph{compatible} if $s_1$ and $s_2$ coincide on $\partial K_1 \cap( k_1 \times[0,1])=\partial K_2 \cap( k_2 \times[0,1])$ for the identification given by the gluing, and if $s$ coincides with the restrictions of $s_1$ and $s_2$ on $\partial K \subset \partial K_1 \cup \partial K_2$.

The main result of the paper is the following. Recall that, for an ideal triangulation of the surface $S$,  $\ZZ_\lambda$ is the square-root Chekhov-Fock algebra defined in \S \ref{sect:CFsquareRoot}.

\begin{figure}[htbp]
\SetLabels
(  .19 * -.2 )  (a) \\
( .81  * -.2 )  (b) \\
( .08 * .63 ) $\epsilon_1$ \\
( .3 * .63 ) $\epsilon_2$ \\
( .73 * .78 )  $\epsilon_1$\\
( .65 * .38 ) $\epsilon_2$ \\
\endSetLabels
\centerline{\AffixLabels{\includegraphics{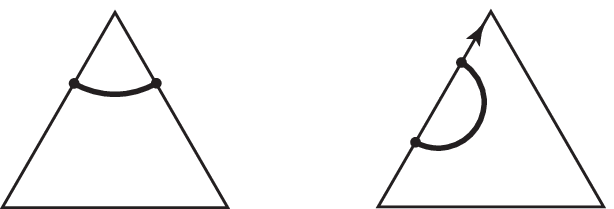}}}
\vskip 15pt
\caption{Elementary skeins in the triangle}
\label{fig:TriangleOneArc}
\end{figure}

\begin{thm}
\label{thm:MainThm}
  For $A=\omega^{-2}$, there  is a unique family of algebra homomorphisms 
$$\Tr_S \col \SSSS(S) \to \ZZ_\lambda,$$
 defined for each surface $S$ and each ideal triangulation $\lambda$ of $S$, such that:
\begin{enumerate}
\item {\upshape (State Sum Property)} If the surface $S$ is obtained by gluing $S_1$ to $S_2$, if  the ideal triangulation $\lambda$ of $S$ is obtained by combining the ideal triangulations $\lambda_1$ of $S_1$ and $\lambda_2$ of $S_2$, and if the skeins $[K_1] \in \SSS(S_1)$ and $[K_2]\in\SSS(S_2)$ are glued together to give $[K] \in \SSS(S)$, then
$$
\quad \quad \quad
 \Tr_S ([K,s]) = \sum_{\text{compatible }s_1, s_2}\Tr_{S_1}([K_1, s_1]) \otimes \Tr_{S_2}([K_2, s_2]) $$
where the sum is over all states $s_1\col \partial K_1 \to \{+,-\}$ and $s_2\col \partial K_2 \to \{+,-\}$ that are compatible with $s\col \partial K \to \{+,-\}$ and with each other. Similarly if the surface $S$, the ideal triangulation $\lambda$ of $S$, and the skein $[K] \in \SSS(S)$ are obtained by gluing  the surface $S_1$, the ideal triangulation $\lambda_1$ of $S_1$, and the skein $[K_1] \in \SSS(S_1)$, respectively,  to themselves, then
$$
\quad \quad
 \Tr_S ([K,s]) = \sum_{\text{compatible }s_1}\Tr_{S_1}([K_1, s_1])  .$$

\item {\upshape (Elementary Cases)} when $S$ is a triangle and $K$ projects to a single arc embedded in $S$, with vertical framing, then 
\begin{enumerate}
\item  in the case of Figure~{\upshape\ref{fig:TriangleOneArc}(a)}, where $\epsilon_1$, $\epsilon_2=\pm$ are the signs associated by the state $s$ to the end points of $K$, then 
$$
\quad\quad \quad \quad \quad
\Tr_S ([K, s]) = 
\left\{  
\begin{aligned}
&0 \text { if } \epsilon_1 =- \text{ and } \epsilon_2 = +\\
& [ Z_1^{\epsilon_1} Z_2^{\epsilon_2} ] \text{ if } \epsilon_1 \not=- \text{ or } \epsilon_2 \not= +
\end{aligned} 
\right.
$$
where $Z_1$ and $Z_2$ are the generators of $\ZZ_\lambda$ associated to the sides $\lambda_1$ and $\lambda_2$ of $S$ indicated, and where $[ Z_1^{\epsilon_1} Z_2^{\epsilon_2} ]=\omega ^{-\epsilon_1\epsilon_2} Z_1^{\epsilon_1} Z_2^{\epsilon_2} =  \omega ^{\epsilon_1\epsilon_2} Z_2^{\epsilon_2} Z_1^{\epsilon_1}$ (identifying the sign $\epsilon = \pm$ to the exponent $\epsilon=\pm1$);

\item in the case of Figure~{\upshape\ref{fig:TriangleOneArc}(b)}, where the end point of $K$ marked by $\epsilon_1$ is higher in $\partial S \times [0,1]$ than the point marked by $\epsilon_2$, 
$$
\Tr_S ([K, s]) = 
\left\{  
\begin{aligned}
0 &\text{ if } \epsilon_1 = \epsilon_2\\
 -\omega^{-5} &\text { if } \epsilon_1 =+ \text{ and } \epsilon_2 = -\\
 \omega^{-1} &\text{ if } \epsilon_1 =- \text{ and } \epsilon_2 = +
\end{aligned} 
\right.
$$
\end{enumerate}
\end{enumerate} 
\end{thm}

In particular, Theorem~\ref{thm:MainThm} immediately gives Theorem~\ref{thm:MainThmIntro} of the introduction by restriction to surfaces without boundary.

In the State Sum Condition (1), note that $\ZZ_\lambda$ is contained in $\ZZ_{\lambda_1} \otimes \ZZ_{\lambda_2}$ when the surfaces $S_1$ and $S_2$ are distinct, and in   $\ZZ_{\lambda_1} $ in the case of a self-gluing. 

In Condition~(2a), the formula $[ Z_1^{\epsilon_1} Z_2^{\epsilon_2} ]= \omega ^{-\epsilon_1\epsilon_2} Z_1^{\epsilon_1} Z_2^{\epsilon_2} =  \omega ^{\epsilon_1\epsilon_2} Z_2^{\epsilon_2} Z_1^{\epsilon_1} $ is more natural than one might think at first glance, as it corresponds to the classical {Weyl quantum ordering} for the monomial $Z_1^{\epsilon_1} Z_2^{\epsilon_2}$. More generally, if $Y_1$, $Y_2$, \dots, $Y_k$ are elements of an algebra such that $Y_iY_j = \omega^{2a_{ij}} Y_jY_i$, the \emph{Weyl quantum ordering} of the monomial $Y_1Y_2\dots Y_k$ is the monomial
$$
[Y_1Y_2 \dots Y_k ]= \omega^{-\sum_{i<j} a_{ij}} Y_1Y_2 \dots Y_k. 
$$
The formula is specially designed to be invariant under all permutations of the $Y_i$.

\subsection{Picture conventions}\label{sect:PictConv}
To work more efficiently with framed links and skeins, we need a convenient way to describe and manipulate  them.

 In practice, we will  represent a link  $K\subset S \times [0,1]$ by its projection to $S$, namely by a 1--dimensional manifold $K'$ immersed in $S$ with $K'\cap \partial S = \partial K'$, and whose only singularities are transverse double points in the interior of $S$; in addition, these double points are endowed with over- or under-crossing information, describing which strand of $K$ lies above the other in $S\times [0,1]$ (with the convention that, when oriented from $0$ to $1$, the $[0,1]$ factor points towards the eye of the reader). 
 
 By adding kinks if necessary, we can always arrange that the framing is vertical at every point of $K$, with the framing vector parallel to the $[0,1]$ factor and pointing towards $1$. 
 
 A crucial information encoded in a framed link $K\subset S \times [0,1]$ is that, for a component $k$ of $\partial S$, the points of $(\partial K) \cap (k\times[0,1])$ are ordered by their elevation. This ordering is not altogether easy to describe on a 2--dimensional picture, and we will resort to the following method to specify this orientation. We choose an arbitrary orientation of $k$. We now have two orderings on $(\partial K) \cap (k\times[0,1])$: one is by order of increasing elevations; the other one is given by the orientation of $k$  if we identify each point of $(\partial K) \cap (k\times[0,1])$ to its projection in $k$. After an isotopy of $K$ (which is elevation-preserving near the boundary), we can always arrange that these two orderings coincide, and we will require this condition to hold in all pictures. 
  
\begin{figure}[htbp]

\SetLabels
\E( .2 * .5 )  $\longrightarrow$\\
\E( .8 * .5 )  $\longrightarrow$\\
\endSetLabels
\centerline{\AffixLabels{\includegraphics{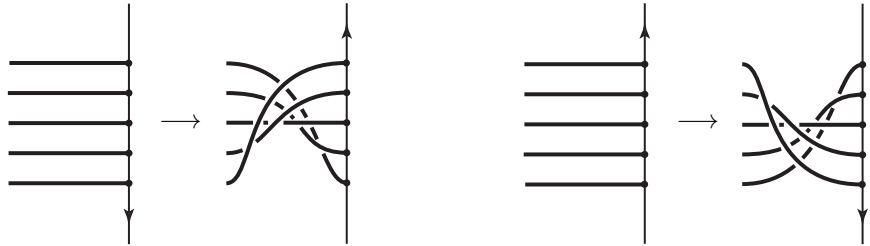}}}

\caption{Reversing a boundary orientation}
\label{fig:HalfTwist}
\end{figure}

 Note that reversing the orientation of $k$ will then oblige us to modify the projection of $K$ by a half-twist near $k$, as in Figure~\ref{fig:HalfTwist}. 
 
With these conventions, the isotopy class of the framed link $K\subset S\times [0,1]$ is then immediately recovered from its projection $K'$ to $S$.

\subsection{Unknots and kinks}
For future reference, we note the following classical facts.

\begin{figure}[htbp]

\SetLabels
\E( .098* .55) \rotatebox{-90}{$\longrightarrow $}\\
\E( .367* .55) \rotatebox{-90}{$\longrightarrow $}\\
\E( .63* .55) \rotatebox{-90}{$\longrightarrow $}\\
\E( .895* .5) \rotatebox{-90}{$\longrightarrow $}\\
( .098*-.3 ) (a)\\
( .367*-.3 ) (b)\\
( .63*-.3 ) (c)\\
( .895*-.3 ) (d)\\
( .17* .68) $K $ \\
( .17* .1) $ K'$ \\
( .435* .68) $K $ \\
( .435* .1) $ K'$ \\
( .7* .68) $K $ \\
( .7* .1) $ K'$ \\
( * ) $ $ \\
\endSetLabels
\centerline{\AffixLabels{ \includegraphics{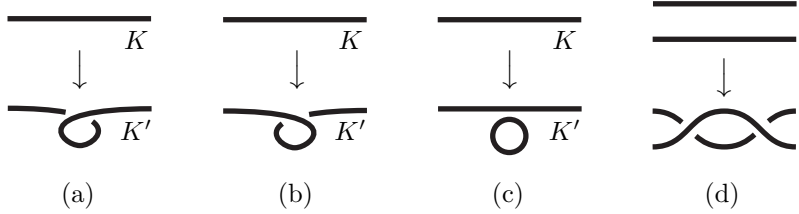} }}
\vskip 20pt
\caption{Adding kinks and unknotted components}
\label{fig:Kinks}
\end{figure}

\begin{lem}
\label{lem:UnknotKink}
If the framed link $K'$ is obtained from $K$ by adding a positive kink as in Figure~{\upshape \ref{fig:Kinks}(a)}, then $[K'] = -A^{-3}[K]$ in $\SSS(S)$. If  $K'$ is obtained from $K$ by adding a negative kink as in Figure~{\upshape \ref{fig:Kinks}(b)}, then $[K'] = -A^{3}[K]$.

If $K'$ is obtained from $K$ by adding a small unknotted circle as in Figure~{\upshape \ref{fig:Kinks}(c)}, then $[K'] = -(A^2+ A^{-2})[K]$ in $\SSS(S)$. 
\end{lem}

In the above statement, we of course assume that those skeins drawn in Figure~\ref{fig:Kinks} follow the picture conventions that we just introduced. 
Adding a positive or negative kink does not change the isotopy class of the link but modifies the framing. 

\begin{proof}
This is an immediate consequence of the skein relations and of the invariance of skeins under the isotopy of Figure~\ref{fig:Kinks}(d). See for instance \cite[Lemmas~3.2 and 3.3]{Licko}. 
\end{proof}

\section{The case of the biangle}
\label{sect:Biangles}

Our proof of Theorem~\ref{thm:MainThm} will make use of \emph{ideal biangles} in addition to ideal triangles. An ideal biangle is the surface $B$ obtained from a closed disk by removing two points from its boundary. In particular, it has two (infinite) edges, and it is also diffeomorphic to the strip delimited by two parallel lines in the plane. 

There is a skein algebra $\SSS(B)$ of the biangle $B$ defined as before. States for skeins are similarly defined.

In this context, we have the following simpler analog of Theorem~\ref{thm:MainThm}.

\begin{figure}[htbp]

\SetLabels
( -.03 * .36 ) $\epsilon_1$ \\
( .22 * .36 ) $\epsilon_2$  \\
( .38 * .57 ) $\epsilon_1$ \\
( .38 * .07 ) $\epsilon_2$  \\
( 1.03 * .57 ) $\epsilon_1$ \\
( 1.03 *  .07) $\epsilon_2$  \\
( .09 * -.3 )  (a) \\
( .5 * -.3 )  (b) \\
( .91 * -.3 ) (c)  \\
( .09 * .8 )  $B$ \\
( .5 * .8 )  $B$ \\
( .91 * .8 ) $B$  \\
\endSetLabels
\centerline{\AffixLabels{\includegraphics{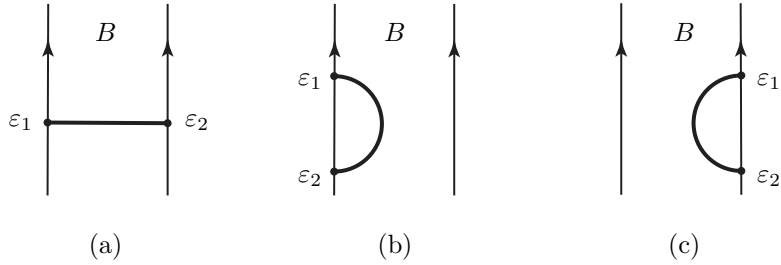}}}
\vskip 20 pt
\caption{Elementary skeins in the biangle}
\label{fig:BiangleOneArc}
\end{figure}

\begin{prop}
\label{prop:BiangleSkein}
Let two numbers $ \alpha, \beta \in \C$ be given, with $ \alpha^2 + \beta^2 = A^5 + A$ and $ \alpha \beta = -A^3$.   Then, there  is a unique family of  algebra homomorphisms 
$$\Tr_B \col \SSSS(B) \to \C,$$
defined for all oriented biangles $B$, such that:  

\begin{enumerate}

\item {\upshape (State Sum Property)} if the biangle $B$ is obtained by gluing together two distinct biangles $B_1$ and $B_2$,  and if $[K_1] \in \SSS(B_1)$ and $[K_2]\in\SSS(B_2)$ are glued together to give $[K] \in \SSS(B)$, then
$$
\quad \quad \quad
 \Tr_B ([K,s]) = \sum_{\text{compatible }s_1, s_2}\Tr_{B_1}([K_1, s_1]) \Tr_{B_2}([K_2, s_2]) ,$$
where the sum is over all states $s_1\col \partial K_1 \to \{+,-\}$ and $s_2\col \partial K_2 \to \{+,-\}$ that are compatible with $s\col \partial K \to \{+,-\}$ and with each other; 

\item {\upshape (Elementary Cases)}  if, using the picture conventions of {\upshape \S\ref{sect:PictConv}},  the biangle $B$ is represented by a vertical strip in the plane as in Figure~{\upshape \ref{fig:BiangleOneArc}} and if $K$  projects to a single arc embedded in $B$, then
\begin{enumerate}
\item  in the case of Figure~{\upshape \ref{fig:BiangleOneArc}(a)}, where $\epsilon_1$, $\epsilon_2=\pm$ are the signs associated by the state $s$ to the end points of $K$,
$$
\quad\quad \quad \quad \quad
\Tr_B([K, s]) = 
\left\{  
\begin{aligned}
&1 \text { if } \epsilon_1 = \epsilon_2 \\
&0 \text{ if } \epsilon_1 \not= \epsilon_2 ;
\end{aligned} 
\right.
$$

\item in the case of Figure~{\upshape \ref{fig:BiangleOneArc}(b)}, 
$$
\Tr_B([K, s]) = 
\left\{  
\begin{aligned}
0 &\text{ if } \epsilon_1 = \epsilon_2\\
 \alpha &\text { if } \epsilon_1 =+ \text{ and } \epsilon_2 = -\\
\beta &\text{ if } \epsilon_1 =- \text{ and } \epsilon_2 = + .
\end{aligned} 
\right.
$$
\end{enumerate}

\end{enumerate}
\end{prop}

Note that the equations $ \alpha^2 + \beta^2 = A^5 + A$ and $ \alpha \beta = -A^3$ only admit four solutions, namely $(\alpha, \beta) = \pm (A^{\frac52}, -A^{\frac12})$ and $ \pm (A^{\frac12}, -A^{\frac52})$. 

\begin{proof} In the case considered, the homomorphism $\Tr_B$ is essentially a version of the Kauffman bracket for tangles. In particular, everything here is fairly classical. However, it is useful to go through the details of the construction to see where the hypotheses on $\alpha$ and $\beta$ come up. 

We will split the proof of Proposition~\ref{prop:BiangleSkein} into several steps. 
We begin with a lemma.

\begin{lem} 
\label{lem:BiangleSkeinStateSum:bd}
For a family of homomorphisms $\Tr_B$ satisfying the properties of Proposition~{\upshape \ref{prop:BiangleSkein}}, then necessarily
$$
\Tr_B([K, s]) = 
\left\{  
\begin{aligned}
0 &\text{ if } \epsilon_1 = \epsilon_2\\
-A^{-3} \alpha &\text { if } \epsilon_1 =+ \text{ and } \epsilon_2 = -\\
-A^{-3} \beta &\text{ if } \epsilon_1 =- \text{ and } \epsilon_2 = +
\end{aligned} 
\right.
$$
when $K$ is as in  Figure~{\upshape \ref{fig:BiangleOneArc}(c)}. 
\end{lem}

\begin{proof}
The proof is provided by Figure~\ref{fig:BiangleRotation}. The equivalence of Figures~\ref{fig:BiangleRotation}(a) and \ref{fig:BiangleRotation}(b) is just obtained by rotating  $B$ by 180 degrees. Reversing the boundary orientations then introduces a half-twist as in Figure~\ref{fig:HalfTwist}, which gives the skein of Figure~\ref{fig:BiangleRotation}(c). Removing the kink, Lemma~\ref{lem:UnknotKink} then shows that this skein is equal to $-A^{-3}$ times the skein of Figure~\ref{fig:BiangleOneArc}(b).  The result then follows from Property~(2b) of Proposition~\ref{prop:BiangleSkein}.
\end{proof}

\begin{figure}[htbp]

\SetLabels
( .22 * .58 ) $\epsilon_1$ \\
( .22 * .08 ) $\epsilon_2$  \\
( .38 * .88 ) $\epsilon_2$ \\
( .38 * .38 ) $\epsilon_1$  \\
( .78 * .58 ) $\epsilon_1$ \\
( .78 *  .08) $\epsilon_2$  \\
\E( .3 * .5) $=$\\
\E( .7 * .5) $=$\\
( .09 * -.3 )  (a) \\
( .5 * -.3 )  (b) \\
( .91 * -.3 ) (c)  \\
\endSetLabels
\centerline{\AffixLabels{\includegraphics{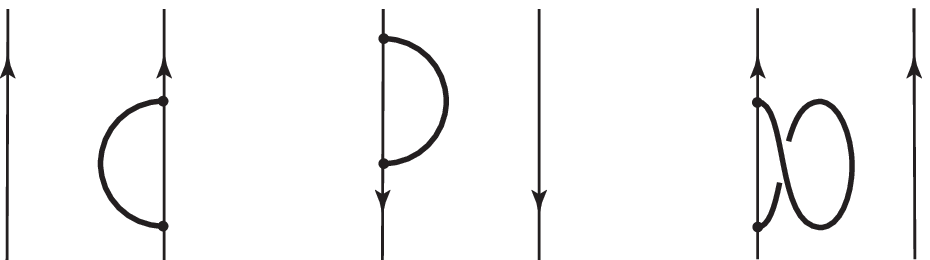}}}
\vskip 20 pt
\caption{The proof of Lemma~\ref{lem:BiangleSkeinStateSum:bd}}
\label{fig:BiangleRotation}
\end{figure}

From now on, when representing a skein in a biangle $B$, we will use the conventions of \S \ref{sect:PictConv} where the two boundary components of $B$ are oriented in a parallel way, as in Figure~\ref{fig:BiangleOneArc}. 

We now prove the uniqueness of the homomorphisms $\Tr_B$.

\begin{lem} 
\label{lem:BiangleSkeinStateUnique}
If there exists a family of homomorphisms $\Tr_B$ satisfying the properties of Proposition~{\upshape\ref{prop:BiangleSkein}} then it is unique.
\end{lem}

\begin{proof}
We first restrict attention to a skein $[K]\in \SSS(B)$ that is represented by a family of arcs and curves \emph{without crossings} in $B$. By general position, isotop $K$ so that it is in bridge position namely so that, as we sweep $B$ from one boundary component to the other, the local maxima and minima are generic and occur at distinct positions.  We can then subdivide $B$ into a union of biangles $B_1$, \dots, $B_n$ so that each $K_i = K \cap (B_i \times [0,1])$ contains at most one maximum or minimum.  Each $K_i$ then is of one of the three types  pictured in Figure \ref{Fig:BiangleSkeinStateSum:Ki}.

\begin{figure}[htbp]

\SetLabels
( .095 * -.2 )  (a) \\
( .5 * -.2 )  (b) \\
( .91 * -.2 ) (c)  \\
(.5 * .71) $\dots$ \\
\R\E(0.005 * .43) $m \left\{\raisebox{.6cm}[.6cm][.6cm]{} \right.$\\
\R\E(.405 * .72) $n \left\{\raisebox{.1cm}[.1cm][.1cm]{} \right.$\\
\R\E(.405 * .17) $m \left\{\raisebox{.1cm}[.1cm][.1cm]{} \right.$\\
\R\E(.81 * .72) $n \left\{\raisebox{.1cm}[.1cm][.1cm]{} \right.$\\
\R\E(.81 * .17) $m \left\{\raisebox{.1cm}[.1cm][.1cm]{} \right.$\\
(.5 * .165) $\dots$ \\
(.91 * .71) $\dots$ \\
(.91 * .165) $\dots$ \\
(.095 * .44) $\dots$ \\
\endSetLabels
\centerline{\AffixLabels{\includegraphics{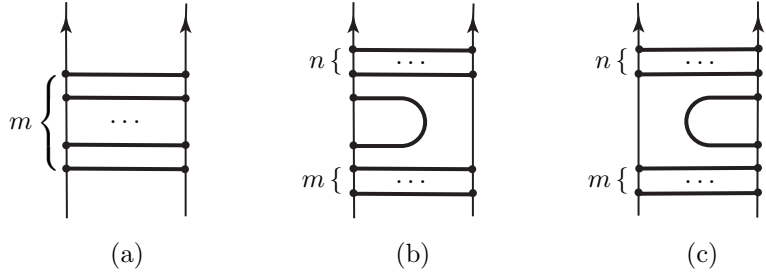}}}
\vskip 20 pt
\caption{The subdivision skeins $K_i$}
\label{Fig:BiangleSkeinStateSum:Ki}
\end{figure}

In particular, if $a$, $b$, $c$ are the three skeins respectively represented in Figures~\ref{fig:BiangleOneArc}(a), \ref{fig:BiangleOneArc}(b) and \ref{fig:BiangleOneArc}(c), $[K_i] = [a]^m$, $[a]^m [b] [a]^n$ or $[a]^m [c] [a]^n$ in the algebra $\SSS(B_i)$. As a consequence, for any state $s_i$, the image of $[K_i, s_i]\in\SSSS (B_i)$ under the algebra homomorphism $\Tr_{B_i}$ is uniquely dedetermined by Properties~(2ab) of Proposition~\ref{prop:BiangleSkein} and by  Lemma \ref{lem:BiangleSkeinStateSum:bd}.

 By the State Sum Property, 
$$
\Tr_B([K,s]) = \sum_{compatible\, s_i} \Tr_{B_1}([K_1, s_1]) \dots \Tr_{B_n}([K_n, s_n]),
$$
so that $\Tr_B([K,s])$ is uniquely determined.

In the case of a skein with crossings, the Kauffman skein relations allow $[K] \in \SSS(B)$ to be expressed as a linear combination  of skeins without crossings.  By linearity of $\Tr_B$ and by uniqueness in the case without crossings, $\Tr_B[K,s]$ is uniquely determined in this general case as well.

This proves Lemma~\ref{lem:BiangleSkeinStateUnique}, namely the uniqueness of the homomorphisms $\Tr_B$. 
\end{proof}

We now  demonstrate the existence of the homomorphisms $\Tr_B$.

First consider the case of a link $K \subset B \times [0,1]$  whose projection to $B$ has no crossing. As in the proof of Lemma~\ref{lem:BiangleSkeinStateUnique}, put $K$ in bridge position, and decompose $B$ as a union of biangles $B_i$ such that $K_i=K \cap (B_i\times[0,1])$ has at most one local maximum or one local minimum for the sweep. 

In this case with no crossing, define
\begin{equation}
\label{eqn:BiangleSkeinStateSum:DefnTr}
\Tr_B([K,s]) = \sum_{compatible\, s_i} \Tr_{B_1}([K_1, s_1]) \dots \Tr_{B_n}([K_n, s_n]),
\end{equation}
where each $\Tr_{B_1}([K_i, s_i]) $ is defined by Conditions~(2ab) of Proposition~\ref{prop:BiangleSkein} or by Lemma~\ref{lem:BiangleSkeinStateSum:bd}.

 Note that, in the above sum, there are few states $s_i$ for which  $\Tr_{B_i}([K_i, s_i])$ is non-zero.  
 
\begin{lem}
\label{lem:BiangleNoCrossing}
For   a skein $[K] \in \SSS(B)$ without crossing, the number $\Tr_B([K,s])$ defined above is  independent of the subdivision of $B$ into biangles $B_i$.
\end{lem}

\begin{proof}
For a given bridge position of $K$, the only freedom in the choice of the biangles $B_i$ is that we can successively add or delete biangles $B_i$ where $K_i$ has no local maximum or minimum. Because of the definition of the corresponding $\Tr_B([K_i,s_i])$ by Condition~(2a) of Proposition~\ref{prop:BiangleSkein},  the state sum  (\ref{eqn:BiangleSkeinStateSum:DefnTr}) providing $\Tr_B([K,s])$ remains unchanged if we add or delete such a $B_i$. 

It therefore suffices to prove independence under the bridge position. By general position, any two bridge positions are related to each other by a sequence of the following moves:
\begin{enumerate}
\item the ``Snake Move'' of Figure~\ref{fig:Snake}, where a local maximum and a local minimum collide and cancel out; this Snake Move actually comes in two types, related to each other by a reflection, according to whether the local maximum sits above or below the local minimum just before the collision; 
\item the inverse of the snake move, which creates a pair of a local maximum and a local minimum;
\item the ``Time Switch Move'', where the times at which two different local extrema occur during the sweep of $B$ are switched; there are 4 types of such Time Switch Moves and their inverses (according to whether they involve local maxima or minima), two of which are represented in Figures~\ref{fig:TimeSwitch1} and \ref{fig:TimeSwitch2}. 
\end{enumerate}

\begin{figure}[htbp]

\SetLabels
\E( .5*.5 ) $ \longleftrightarrow$ \\
\endSetLabels
\centerline{\AffixLabels{ \includegraphics{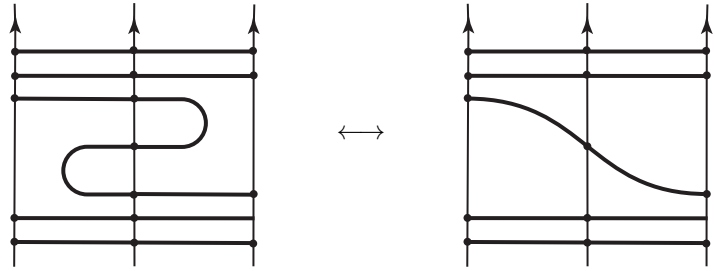} }}
\caption{The Snake Move and its inverse}
\label{fig:Snake}
\end{figure}

\begin{figure}[htbp]

\SetLabels
\E( .5*.5 ) $ \longleftrightarrow$ \\
\endSetLabels
\centerline{\AffixLabels{ \includegraphics{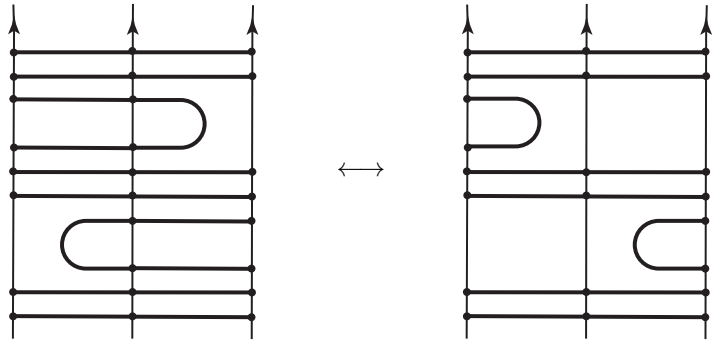} }}
\caption{A Time Switch Move}
\label{fig:TimeSwitch1}
\end{figure}

\begin{figure}[htbp]

\SetLabels
\E( .5*.5 ) $ \longleftrightarrow$ \\
\endSetLabels
\centerline{\AffixLabels{ \includegraphics{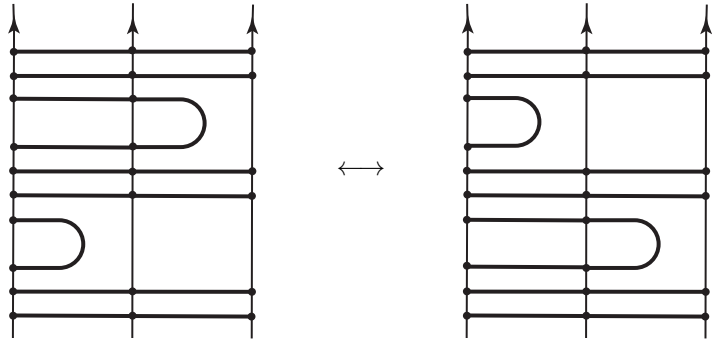} }}
\caption{Another Time Switch Move}
\label{fig:TimeSwitch2}
\end{figure}

The invariance of the state sum formula (\ref{eqn:BiangleSkeinStateSum:DefnTr}) for $\Tr_B([K,s])$ under these Snake and Time Switch Moves follows from an easy computation using the definition of the $\Tr_{B_i}([K_i,s_i])$ by Lemma~\ref{lem:BiangleSkeinStateSum:bd} and Condition~(2) of Proposition~\ref{prop:BiangleSkein}. In particular, the hypothesis that  $\alpha\beta = -A^3$ is critical for the Snake Move. 

This proves Lemma~\ref{lem:BiangleNoCrossing}, and uniquely defines $\Tr_B([K,s])$ for every skein with no crossing. 
\end{proof}

In particular:

\begin{lem} \label{lem:BiangleSkeinStateSum:circle}
If $K$ projects to a small circle embedded in the biangle $B$, $\Tr_B([K, \varnothing])= -A^2 - A^{-2}$.
\end{lem}
\begin{proof} Putting $K$ in bridge position with only one local maximum and one local minimum, the state sum formula  (\ref{eqn:BiangleSkeinStateSum:DefnTr}) involves only four compatible states, two of which contribute 0 to the sum. This gives,
$$
\Tr_B([K, \varnothing] )= 0 +0 + \alpha(-A^{-3} \alpha) + \beta(-A^{-3}\beta) = -A^2 - A^{-2},
$$
using the hypothesis that $\alpha^2 +\beta^2 = A + A^5$.
 \end{proof}

We now define $\Tr_B([K,s])$ for an arbitrary link $K$ with a state $s$. By resolving all the crossings of $K$ and applying the skein relation, write $[K] \in \SSS(B)$ as a linear combination 
$$
[K] = \sum_{i=1}^k A^{n_i} [K_i]
$$
 of skeins $[K_i]$ where the link $K_i$ has no crossing. Then, define
$$
\Tr_B([K,s]) =\sum_{i=1}^k A^{n_i}  \Tr_B([K_i, s_i]).
$$

\begin{lem}
\label{lem:TraceIsotInv}
The number $\Tr_B([K,s])$ defined above is independent of the framed isotopy class of $K$.
\end{lem}
\begin{proof}
It suffices to show invariance under the second and third Reidemeister Moves. This is a classical consequence of Lemma~\ref{lem:BiangleSkeinStateSum:circle} (see for instance \cite[Lem\-ma~3.3]{Licko}). 
\end{proof}

By construction, it is immediate that the $\Tr_B([K,s])$ satisfy the skein relation. Therefore, the construction provides a linear map
$$
\Tr_B \col \SSSS(B) \to \C. 
$$
It is also  immediate that this linear map also satisfies the State Sum Property (1) of Proposition~\ref{prop:BiangleSkein}. It remains to show that it is an algebra homomorphism. 

\begin{lem}
For any two stated skeins $[K, s]$, $[K', s']\in \SSSS(B)$, 
$$
\Tr_B ([K, s][K', s']) = \Tr_B([K,s]) \,\, \Tr_B([K', s']).
$$
\end{lem}

\begin{figure} [htb]

\SetLabels
\E(.25 *.27 ) $ K$ \\
\E( .75* .74) $K' $ \\
\E(.1 *.27 ) $\dots $ \\
\E(.42 *.27 ) $\dots $ \\
\E(.75 *.27 ) $\dots $ \\
\E(.9 *.73 ) $\dots $ \\
\E(.58 *.73 ) $\dots $ \\
\E(.25 *.73 ) $\dots $ \\
\endSetLabels
\centerline{\AffixLabels{\includegraphics{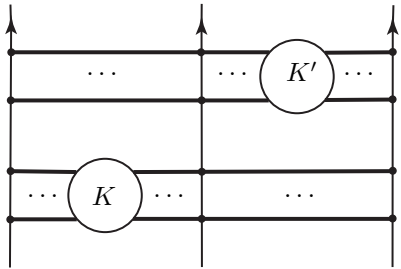}}}

\caption{Superposition of skeins $K$ and $K'$}

 \label{Fig:BiangleSkeinStateSum:KK'}
\end{figure}

\begin{proof}
Our convention of using parallel orientations for the two boundary components of $B$ turns out to be convenient here.
The product $[K, s][K', s']$ is equal to $[KK', s \cup s']$, where $KK'$ denotes the superposition of the links $K$ and $K'$. 
Because of the orientation convention, the superposition $KK'$ of $K$ and $K'$ can be isotoped so that $K$ and $K'$ sit side by side in $B$, with $K'$ above $K$ on the sheet of paper as in Figure~\ref{Fig:BiangleSkeinStateSum:KK'}. If we use this configuration in our construction of $\Tr_B([KK', s \cup s']) $, it is then immediate that $\Tr_B([KK', s \cup s']) = \Tr_B([K,s]) \, \Tr_B([K', s']) $.
\end{proof}

This completes the proof of Proposition~\ref{prop:BiangleSkein}. 
\end{proof}

Although the definition of $\Tr_B$ may seem complicated,  its computation is much simpler in practice.  Indeed, if $K$ is a link whose projection to $B$ has no crossing, each of its components is a closed curve, or an arc of one of the three types of Figure~\ref{fig:BiangleOneArc}.  If, in addition, $K$ is endowed with a state $s$ and if $\epsilon_1=\pm$ and $\epsilon_2=\pm$, let $a^{\epsilon_1}_{\epsilon_2}$ be the number of components of the type of Figure~\ref{fig:BiangleOneArc}(a) where the state $s$ assigns signs  $\epsilon_1$ and $\epsilon_2$ as in that figure; let $b^{\epsilon_1}_{\epsilon_2}$ be the number of components of the type of Figure~\ref{fig:BiangleOneArc}(b); let $c^{\epsilon_1}_{\epsilon_2}$ be the number of components of the type of Figure~\ref{fig:BiangleOneArc}(c); and let $d$ be the number of closed components of $K$. 

\begin{lem}
\label{lem:TraceSimpleDef}
For a stated skein $[K,s] \in \SSSS(B) $ with no crossing, let $a^{\epsilon_1}_{\epsilon_2}$, $b^{\epsilon_1}_{\epsilon_2}$, $c^{\epsilon_1}_{\epsilon_2}$ and $d$ be defined as above. If one of $a^+_-$, $a^-_+$, $b^+_+$, $b^-_-$, $c^+_+$, $c^-_-$ is non-zero, then $\Tr_B([K,s]) = 0$. Otherwise, 
$$
\Tr_B([K,s]) = \alpha{\vrule height 1.8ex width 0pt}^{b^+_-}\,\beta{\vrule height 2ex width 0pt}^{b^-_+} \bigl(-A^{-3} \alpha \bigr)^{c^+_-} \bigl(-A^{-3}\beta\bigr)^{c^-_+} \bigl(-A^2-A^{-2} \bigr)^d.
$$
\end{lem}
\begin{proof}
Isotop $K$ so that it is in bridge position and so that: each arc component of $K$ has only one local maximum or minimum; the projection of each closed component of $K$ to $B$ bounds a disk whose interior is disjoint from the projection of $K$. The formula then follows from the definition of $\Tr_B([K,s])$. 
\end{proof}

In particular, for a skein with no crossings, $\Tr_B([K,s]) $ is independent of the relative nesting of the components of the projection of $K$ to $B$. 

The following two observations will be useful later on.

\begin{lem}
\label{lem:NonZeroBracket}
If $\Tr_B([K,s]) $ is non-zero, the stated skein $[K,s] \in \SSSS(B) $ is \emph{balanced} is the sense that the sum of the signs assigned by $s$ to the components of $\partial K$ in one component of $\partial B$ is equal to the sum of the signs in the other component of $\partial B$. 
\end{lem}
\begin{proof}
Lemma~\ref{lem:TraceSimpleDef} proves this for skeins with no crossings. The general case follows from this one by resolving all the crossings and applying the skein relations. 
\end{proof}

\begin{figure}[htb]
\SetLabels
\R( 0* .72) $\epsilon_1 $ \\
\R( 0* .24) $\epsilon_2 $ \\
\L( 1* .72) $\epsilon_1 '$ \\
\L( 1* .24) $\epsilon_2 '$ \\
\endSetLabels
\centerline{\AffixLabels{ \includegraphics{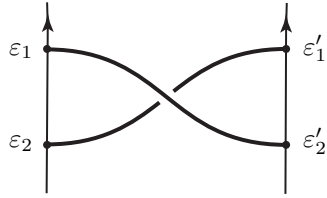} }}

\caption{A right-handed crossing}
\label{Fig:BiangleSkein:RightTwist}
\end{figure}

The following is an immediate computation, based on the definitions. 
\begin{lem}
\label{lem:RightTwist}
For the right half-twist $[K,s]$ represented in Figure~{\upshape\ref{Fig:BiangleSkein:RightTwist}}, 
$$
 \Tr_B[K,s]  = 
\begin{cases} 
A & \text{if  } \epsilon_1 = \epsilon'_1=\epsilon_2 = \epsilon'_2\\
A- A^{-4} \alpha^2& \text{if  }  \epsilon_1 = \epsilon'_1=+  \text{ and } \epsilon_2 = \epsilon'_2=- \\
A- A^{-4} \beta^2& \text{if  }  \epsilon_1 = \epsilon'_1=- \text{ and }  \epsilon_2 = \epsilon'_2 =+ \\
-A^{-4}\alpha\beta& \text{if  } \epsilon_1 = \epsilon'_2 \neq \epsilon_2 = \epsilon'_1\\
0 &\text{otherwise.}
\end{cases}
$$\qed
\end{lem}

\section{Split ideal triangulations}

A \emph{split ideal triangulation} $\widehat \lambda$ is obtained from an ideal triangulation $\lambda$ by replacing each edge of $\lambda$ by two parallel copies of it, separated by a biangle. In particular, $\widehat \lambda$ is a cell decomposition of $S$ whose faces consists of finitely many triangles $T_1$, $T_2$, \dots, $T_m$ (each corresponding to a face of $\lambda$) and finitely many biangles $B_1$, $B_2$, \dots, $B_n$ (each corresponding to an edge of $\lambda$). 

\begin{lem}
\label{lem:GoodPosition}
Let $K$ be a framed link in $S \times [0,1]$ and let $\widehat\lambda$ be a split ideal triangulation of $S$. Then $K$ can be isotoped so that:
\begin{enumerate}
\item for every edge $\widehat\lambda_i$ of $\widehat\lambda$, $K$ is transverse to $\widehat\lambda_i \times [0,1]$;
\item for every triangle face $T_j$ of $\widehat \lambda$, $K \cap(T_j \times [0,1])$ consists of finitely many disjoint arcs, each of which is contained in a constant elevation surface $S \times *$ and joins two distinct components of $\partial T_i \times [0,1]$;
\item for every triangle face $T_j$ of $\widehat \lambda$, the components of $K \cap(T_j \times [0,1])$ lie at different elevations, and their framings are vertical pointing upwards.
\end{enumerate}
\end{lem}

The effect of Lemma~\ref{lem:GoodPosition} is to push all the complexities of $K$ into the part of $S\times[0,1]$ that lies above the biangles of $\widehat \lambda$. 

\begin{proof}
Select a spine $Y_j$ for each ideal triangle $T_j$, namely an infinite Y-shaped subset such that $T_j$ properly collapses on $Y_j$, as in Figure~\ref{fig:TriangleSpine}.

\begin{figure}[htb]

\SetLabels
( .5 * .05 )  $T_j$ \\
( .55 * .4 ) $Y_j$ \\
\endSetLabels
\centerline{\AffixLabels{\includegraphics{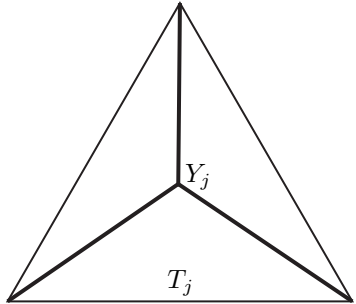}}}

\caption{The spine of an ideal triangle}
\label{fig:TriangleSpine}
\end{figure}

By generic position we can arrange that, for each spine $Y_j$, the link $K$ is disjoint from the singular locus $\{v_j\}  \times[0,1]$ of each $Y_j \times[0,1]$, corresponding to the 3--valent vertex $v_j \in Y_j$, and transverse to the rest of $Y_j \times[0,1]$. With a further isotopy we can assume that on a small neighborhood  $U_j$ of $Y_j \times[0,1]$ each component of $K\cap U_j$ has constant elevation, and that distinct components have distinct elevations. Finally, the framing can be modified so that it is vertical and pointing upwards on $K\cap U_j$. 

By definition of the spines $Y_j$, the union of the $T_j\times[0,1]$ can be isotoped inside the union of the $U_j$, and this by an isotopy  which respects all level surfaces $S\times *$ and which sends vertical arc $*\times[0,1]$ to vertical arc. Modifying $K$ by the inverse of this isotopy puts it in the required position. 
\end{proof}

When $K$ satisfies the conclusions of Lemma~\ref{lem:GoodPosition}, we will say that it is in \emph{good position} with respect to the split ideal triangulation $\widehat\lambda$. 

\begin{figure}[htbp]

\includegraphics{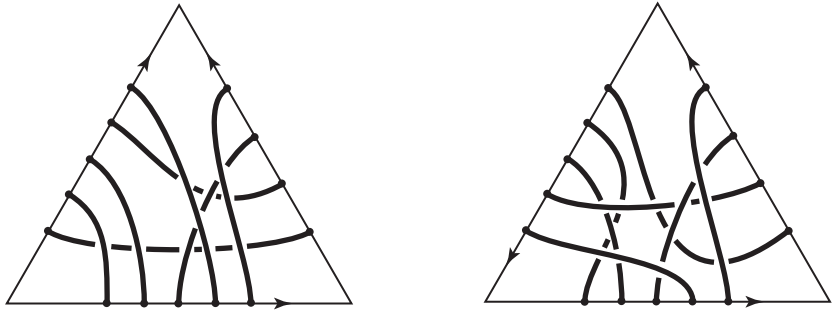}
\caption{A link in good position over a triangle}
\label{fig:TriangleExample}
\end{figure}

When $K$ is in good position, we can always isotop it so that, over each triangle face $T_j$, $K\cap (T_j \times [0,1])$ projects to finitely many disjoint arcs embedded in $T_j$. However, the projection is usually  much more complicated when we use the conventions of \S \ref{sect:PictConv} to represent the ordering of the components of $K\cap (T_j \times [0,1])$. The two pictures of Figure~\ref{fig:TriangleExample} illustrate the same example of a link in good position over a triangle $T_j$, drawn with the conventions of \S \ref{sect:PictConv} but with two different choices of orientations for the sides of $T_j$.

\begin{figure}[htbp]
\SetLabels
\E( .5 *  .5) $\longleftrightarrow$ \\
(  *  )  \\
(  *  )  \\
(  *  )  \\
\endSetLabels
\centerline{\AffixLabels{\includegraphics{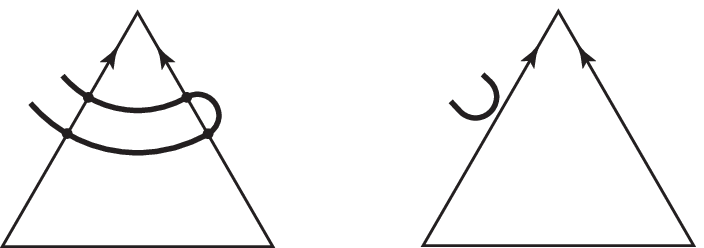}}}
\caption{Move (I)}
\label{fig:move1}
\end{figure}

\begin{figure}[htbp]
\SetLabels
\E( .5 *  .5) $\longleftrightarrow$ \\
(  *  )  \\
(  *  )  \\
(  *  )  \\
\endSetLabels
\centerline{\AffixLabels{\includegraphics{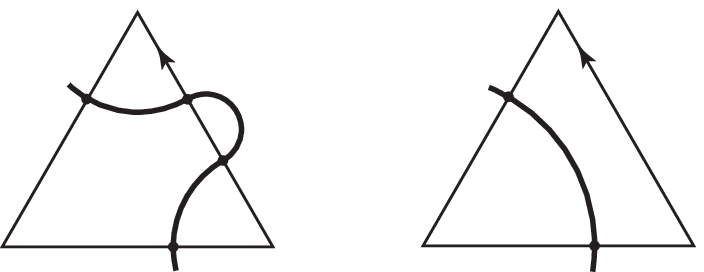}}}
\caption{Move (II)}
\label{fig:move2}
\end{figure}

\begin{figure}[htbp]
\SetLabels
\E( .5 *  .5) $\longleftrightarrow$ \\
(  *  )  \\
(  *  )  \\
(  *  )  \\
\endSetLabels
\centerline{\AffixLabels{\includegraphics{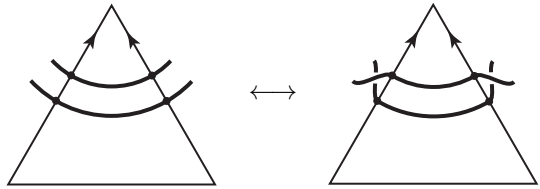}}}
\caption{Move (III)}
\label{fig:move3}
\end{figure}

\begin{figure}[htbp]
\SetLabels
\E( .5 *  .5) $\longleftrightarrow$ \\
(  *  )  \\
(  *  )  \\
(  *  )  \\
\endSetLabels
\centerline{\AffixLabels{\includegraphics{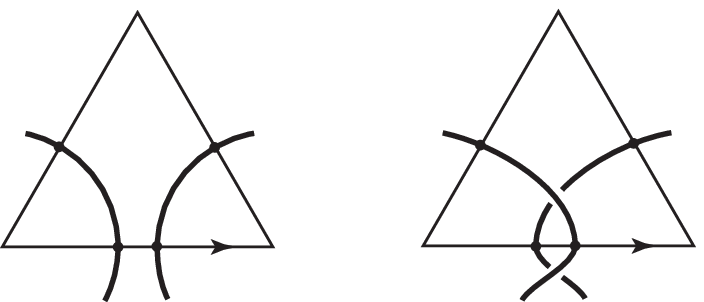}}}
\caption{Move (IV)}
\label{fig:move4}
\end{figure}

\begin{figure}[htbp]
\SetLabels
\E( .5 *  .5) $\longleftrightarrow$ \\
(  *  )  \\
(  *  )  \\
(  *  )  \\
\endSetLabels
\centerline{\AffixLabels{\includegraphics{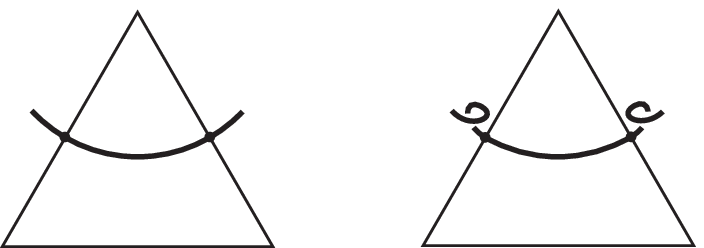}}}
\caption{Move (V)}
\label{fig:move5}
\end{figure}

Figures~\ref{fig:move1}-\ref{fig:move5} describe five moves, occurring in a neighborhood of a triangle $T_j$, which isotop $K$ from one good position to another. It is understood there that $K$ is in good position with respect to the split ideal triangulation, and in particular that each component of $K \cap (T_j \times[0,1])$ has constant elevation. The intersection $K \cap (T_j \times[0,1])$ may include more components than the zero, one or two arcs shown. However, when two arcs are shown, it is understood that they are adjacent to each other for the ordering of the components of $K \cap (T_j \times[0,1])$ given by their elevation; their ordering with respect to each other is determined by the ordering of their end points, indicated by the arrows with the conventions of \S \ref{sect:PictConv}.

Moves (I) and (II) of Figures~\ref{fig:move1}-\ref{fig:move2} eliminate a U-turn in  biangles adjacent to $T_j$. Moves (III) and (IV) transpose the elevations of two components of $K \cap (T_j \times[0,1])$. The kinks added by Move (V) change the framing in two biangles adjacent to $T_j$ (if these two biangles are distinct), but not the isotopy class of the global framing of $K$ since the two kinks have opposite signs.

\begin{lem}
\label{lem:GoodPositionMoves}
Let $K$ and $L$ be two framed links in $S\times[0,1]$ satisfying the conclusions of Lemma~{\upshape\ref{lem:GoodPosition}}, namely in good position with respect to $\widehat\lambda$. If $K$ and $L$ are isotopic, then they can be connected by a sequence of framed links $K=K_0$, $K_1$, $K_2$, \dots, $K_l$, $K_{l+1}=L$ such that each $K_i$ is in good position with respect to $\widehat\lambda$, and such that one goes from $K_i$ to $K_{i+1}$ by either an isotopy keeping the link in good position with respect to $\widehat\lambda$, or one of the moves {\upshape (I)--(V)} of Figures~{\upshape\ref{fig:move1}--\ref{fig:move5}} and their inverses.  
\end{lem}

\begin{proof} Choose a spine $Y_j$ in each triangle face $T_j$ of the split ideal triangulation $\widehat\lambda$, as in the proof of Lemma~\ref{lem:GoodPosition}. For a framed link $K$ in $S\times [0,1]$, this proof of Lemma~\ref{lem:GoodPosition} shows that being in good position with respect to the split triangulation $\widehat \lambda$ is essentially equivalent to the property that, for each spine $Y_j$, 
\begin{enumerate}
\item the link $K$ is transverse to $Y_j \times [0,1]$, and in particular is disjoint from the singular locus $\{v_j\}  \times [0,1]$, where $v_j$ is the trivalent vertex of $Y_j$;
\item the elevations of the points of $K \cap (Y_j \times[0,1])$ are distinct;
\item the framing of $K$ is vertical, pointing upwards, at each point of $K \cap (Y_j \times[0,1])$.
\end{enumerate}
When these properties hold, we say that $K$ is \emph{in good position with respect to the spine $Y_j$}. 

Indeed, by adjusting $K$ in a neighborhood of the $Y_j$ as in the proof of Lemma~\ref{lem:GoodPosition} and by collapsing the triangles $T_j$ in this neighborhood by isotopies, one easily goes back and forth between links that are in good position with respect to the split triangulation $\widehat\lambda$ and links that are in good position with respect to the spines $Y_j$. 

Let us first neglect the elevations of the intersection points and the framings.

If the two framed links $K$ and $L$ are in good position with respect to the spines $Y_j$, and if they are isotopic, let us choose the isotopy to be generic with respect to the $Y_j \times[0,1]$. Then, during the isotopy,  the link remains transverse to the $Y_j \times[0,1]$ except at finitely many times where, either the link crosses one of the singular loci $\{v_j\} \times[0,1]$, or it is tangent to one of the smooth parts $(Y_j - \{v_j\}) \times[0,1]$. 

When the link crosses the singular locus $\{v_j\} \times[0,1]$, an intersection point of the link with a component of $(Y_j - \{v_j\}) \times[0,1]$ gets replaced with two intersection points, one in each of the two other components of $(Y_j - \{v_j\}) \times[0,1]$, or the converse holds. Translating this in terms of links in good position with respect to the split ideal triangulation $\widehat\lambda$ gives Move~(II) and its inverse. 

When the link becomes tangent to $(Y_j - \{v_j\}) \times[0,1]$ then, generically, either two intersection points of the link with $(Y_j - \{v_j\}) \times[0,1]$ collide and cancel out, or a pair of intersection points gets created in the inverse process. This is described by Move~(I) and its inverse. 

So far, we had neglected the elevation of the intersection points. At finitely many times during the generic isotopy, the elevations of two points of some $K \cap (Y_j \times[0,1])$ will cross each other. This is described by Moves~(III) and (IV), according to whether the points are in the same component of $(Y_j - \{v_j\}) \times[0,1]$ or not. 

In particular, when we care about elevations, Moves~(III) and (IV) enable us to avoid having to consider two versions of Moves~(I) and (II), one for each ordering of the two points of $K \cap (Y_j \times[0,1])$. 

Finally, we have to worry about framings. At some time in the isotopy, we will need to move the framing at a point of $K \cap (Y_j \times[0,1])$ from vertical position to vertical position by rotating it by a certain number of full turns. This is accomplished by several applications of Move~(V) or its inverse. 
\end{proof}

\section{The quantum trace as a state sum}

We now begin the proof of Theorem~\ref{thm:MainThm}.

Let $K$ be a framed link in $S\times[0,1]$, with a state $s\col \partial K \to \{+, -\}$. 
Let $\lambda$ be an ideal triangulation for $S$. 

Let $\widehat\lambda$ be a split ideal triangulation associated to $\lambda$. By an isotopy, put $K$ in good position with respect to $\widehat\lambda$ as in Lemma~\ref{lem:GoodPosition}. The conclusions of this lemma guarantee that, for every triangle face $T_j$ or biangle face  $B_i$ of $\widehat\lambda$, the intersections $K_j=K\cap(T_j\times[0,1])$ and $L_i=K\cap(B_i\times[0,1])$ are framed links in $T_j \times[0,1]$ or $B_i\times[0,1]$. 

Suppose that, in addition to the state $s\col \partial K \to \{+,-\}$, we are given a state $s_j\col \partial K_j \to \{ +,-\}$ for each $K_j$ with $j=1$, $2$, \dots, $m$, and a state $t_i\col \partial L_i \to \{ +,-\}$ for each $L_i$ with $i=1$, $2$, \dots, $n$. Note that exactly two of these states are defined at every point of $\partial K \cup \bigcup_{j=1}^m \partial K_j \cup \bigcup_{i=1}^n \partial L_i$. We say that all these states $s$, $s_j$ and $t_i$ are \emph{compatible} if they coincide whenever they are defined at the same point.  

For a triangle $T_j$, let $k_1$, $k_2$, \dots, $k_l$ be the components of $K_j$, in order of increasing elevation (remember that the elevation is constant on each $k_i$, and that distinct $k_i$ have distinct elevations). Then $K_j=k_1k_2\dots k_l$ in the link algebra $\mathcal K(T_j)$; note that the order of the terms in this product is important. Let $\Tr_{T_j} (k_i, s_j) \in \ZZ_{T_j}$ be defined as in (2a) of Theorem~\ref{thm:MainThm}. Then, define
$$
\Tr_{T_j} (K_j, s_j) = \Tr_{T_j} (k_1, s_j)\Tr_{T_j} (k_2, s_j)\dots \Tr_{T_j} (k_l, s_j)  \in \ZZ_{T_j}.
$$

For a biangle $B_i$ of $\widehat \lambda$, let $\Tr_{B_i} (L_i, t_i) \in \C$ be the scalar provided by Proposition~\ref{prop:BiangleSkein}. 

We can then consider the tensor product
$$
\prod_{i=1}^n \Tr_{B_i} (L_i, t_i) \, \bigotimes_{j=1}^m\Tr_{T_j} (K_j, s_j) \in \bigotimes_{j=1}^m  \ZZ_{T_j}
$$

Recall that the triangles $T_j$ are identified to the faces of the ideal triangulation $\lambda$, and that the Chekhov-Fock square root algebra $\ZZ_\lambda$ is also contained in $\bigotimes_{j=1}^m  \ZZ_{T_j}$. 

\begin{lem}
\label{lem:QuantumTraceInCheFock}
If the states $s$, $s_j$ and $t_i$ are compatible, the element 
$$
\prod_{i=1}^n \Tr_{B_i} (L_i, t_i)\,  \bigotimes_{j=1}^m\Tr_{T_j} (K_j, s_j) \in \bigotimes_{j=1}^m  \ZZ_{T_j}
$$
 is contained in the square root Chekhov-Fock algebra $\ZZ_\lambda$ of {\upshape \S\ref{sect:CFsquareRoot}}. 
\end{lem}

\begin{proof}
We first have to check that, when the monomial
$$
\prod_{i=1}^n \Tr_{B_i} (L_i, t_i) \, \bigotimes_{j=1}^m\Tr_{T_j} (K_j, s_j)
$$
is different from $0$, the generators $Z_{ja}$ and $Z_{ka}$ associated to the two sides of an edge $\lambda_i$ of $\lambda$ appear with the same exponent in the monomial. This is an immediate consequence of Lemma~\ref{lem:NonZeroBracket} and of the definition of the terms $\Tr_{T_j} (K_j, s_j)$. 

The fact that this monomial satisfies the parity condition defining the square root algebra $\ZZ_\lambda$ automatically follows from the definitions.
\end{proof}

Define
$$
\Tr_S (K, s) = \sum_{\text{compatible }s_j, t_i} \ 
\prod_{i=1}^n \Tr_{B_i} (L_i, t_i)  \bigotimes_{j=1}^m\Tr_{T_j}  (K_j, s_j) \in \ZZ_\lambda(K),
$$
where the sum of over all choices of states $s_j\col \partial K_j \to \{ +,-\}$, $j=1$, $2$, \dots, $m$, and  $t_i\col \partial L_i \to \{ +,-\}$, $i=1$, $2$, \dots, $n$, that are compatible with $s\col \partial K \to \{+,-\}$ and with each other. 

The key step in the proof of Theorem~\ref{thm:MainThm} is the following. 

Recall that our definition of $\Tr_B$ for biangles $B$ in \S \ref{sect:Biangles} depended on two parameters $\alpha$ and $\beta$ such that $\alpha^2 + \beta^2=A + A^5$ and $\alpha\beta = -A^3$. This allowed four possibilities $(\alpha, \beta) = \pm (A^{\frac12}, -A^{\frac52})$ and $\pm (A^{\frac52}, -A^{\frac12})$ for these parameters. They now need to be even more restricted. Going over the proof of Proposition~\ref{prop:QuantumTraceWellDefined}, the reader will readily check that these restrictions on $A$, $\alpha$ and $\beta$ are necessary for the statement to hold. 

\begin{prop}
\label{prop:QuantumTraceWellDefined}
If $A=\omega^{-2}$, $\alpha = -\omega^{-5}$ and $\beta = \omega^{-1}$, the
above element 
$$
\Tr_S (K, s)  \in \ZZ_\lambda(K),
$$
depends only on the isotopy class of $K$ and on the state $s$.
\end{prop}

\begin{proof} By Proposition~\ref{prop:BiangleSkein}, $\Tr_S (K, s)$ is invariant under isotopy respecting good position with respect to the split ideal triangulation $\widehat\lambda$. Therefore, we only need to check that it remains unchanged under the Moves (I)--(V) of  Lemma~\ref{lem:GoodPositionMoves}.

These moves involve a triangle $T_j$, adjacent to three biangles $B_{i_1}$, $B_{i_2}$, $B_{i_3}$. We will restrict attention to the case where these three biangles are distinct. Since it involves only minor modifications in notation and no new arguments, we leave as  an exercise to the reader the task of adapting our proof to the case where two of the biangles coincide. 

To alleviate the notation, we can assume that the triangle involved is the triangle $T_1$, while the adjacent biangles are $B_1$, $B_2$, $B_3$. In addition, in each of Figures~\ref{fig:move1}--\ref{fig:move5}, we will assume that the $B_i$ are indexed as in Figure~\ref{fig:LabelConventions}. In particular, the square root  algebra $\ZZ_{T_1}$ is defined by generators $Z_{11}$, $Z_{12}$ and $Z_{13}$, respectively associated to the edges $T_1\cap B_1$, $T_1\cap B_2$ and $T_1\cap B_3$, such that $Z_{1i}Z_{1(i+1)} = \omega^2 Z_{1(i+1)}Z_{1i}$ if we count indices modulo 3.

\begin{figure}[htbp]

\SetLabels
\E( .5* .35) $T_1 $ \\
(.1 *.4 ) $ B_1$ \\
( .9*.4 ) $B_2 $ \\
(.5 * -.15) $ B_3$ \\
\endSetLabels
\centerline{\AffixLabels{\includegraphics{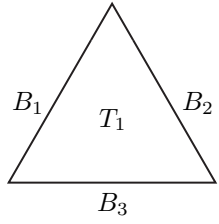}}}
\hskip 10pt
\caption{Labeling conventions}
\label{fig:LabelConventions}
\end{figure}

Consider a move of type (I)--(V), going from a framed link $K$ (on the left of each of Figures~\ref{fig:move1}--\ref{fig:move5}) to a framed link $K'$ (on the right of Figures~\ref{fig:move1}--\ref{fig:move5}).  In the above state sums for $\Tr_S(K, s)$ and $\Tr_S(K', s)$ we group terms so that, in each group, the families of compatible states $s_j$, $t_i$ for $K$  and $s_j'$, $t_i'$ for $K'$ coincide outside of the parts of $K$ and $K'$ shown. We then have to show that, in each group, the sum of the contributions to $\Tr_S(K, s)$ of the compatible states $s_j$, $t_i$ for $K$ considered is equal to the sum of the contributions to $\Tr_S(K', s)$  of the corresponding compatible states $s_j'$, $t_i'$ for $K'$. We will group terms even further according to the powers of the generators $Z_{11}$, $Z_{12}$, $Z_{13}$ of $\ZZ_{T_1}$ (associated to the sides of the triangle $T_1$) contributed by these states.

\begin{figure}[htbp]
\SetLabels
(.465 * .66 ) $+$\\
(.445 * .37 ) $+$\\
(.535 * .66 ) $-$\\
(.555 * .37 ) $+$\\
(.255 * .66 ) $+$\\
(.235 * .37 ) $-$\\
(.325 * .66 ) $+$\\
(.345 * .37 ) $-$\\
(.045 * .66 ) $+$\\
(.025 * .37 ) $+$\\
(.115 * .66 ) $+$\\
(.135 * .37 ) $-$\\
(.675 * .66 ) $-$\\
(.655 * .37 ) $+$\\
(.745 * .66 ) $-$\\
(.765 * .37 ) $+$\\
\E( .19* .5) $+$\\
\E(.40 * .5) $+$\\
\E(  .605* .5) $+$\\
\E( .81* .7) $?$\\
\E(.81 * .5) $=$\\
\endSetLabels
\centerline{\AffixLabels{ \includegraphics{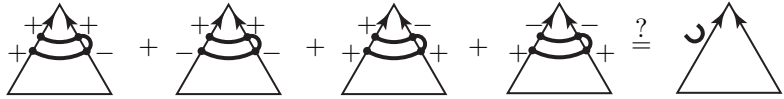} }}

\caption{States for Move (I)}
\label{fig:Move1bis}
\end{figure}

In the case of Move~(I), for a given family of compatible states $s_j'$, $t_i'$ for $K'$, there are $2^4$ families of compatible states $s_j$, $t_i$  for $K$ that coincide with $s_j'$, $t_i'$ outside of the area shown,  but only four of these give a non-zero contribution to $\Tr_S(K,s)$. These are indicated in Figure~\ref{fig:Move1bis}. The equality signs in this figure means that we have to show that the contributions of the terms on one side of the equation add up to the contributions of the other side. 

For the labeling conventions of Figure~\ref{fig:LabelConventions} and remembering that $L_i$ denotes the portion of $K$ that is above the biangle $B_i$, let $t_1^{++}$, $t_1^{-+}$, $t_1^{+-}$ be the states for $L_1$ described by the first, second and fourth triangles in Figure~\ref{fig:Move1bis}, in this order. (To explain the notation, note the signs assigned by these states to the two points of $K\cap B_1 \cap T_1$ shown, for the orientation of the edge $B_1 \cap T_1$ specified by the arrow.) Similarly, let $t_2^{-+}$ and $t_2^{+-}$ be the states for $L_2$ represented in the first and third triangles, respectively.  

By the State Sum Property of Proposition~\ref{prop:BiangleSkein}, 
\begin{align*}
\Tr_{B_1}(L_1', t_1') &= \alpha \,\Tr_{B_1}(L_1, t_1^{-+}) + \beta \, \Tr_{B_1}(L_1, t_1^{+-})\\
&= - \omega^{-5} \,\Tr_{B_1}(L_1, t_1^{-+}) +  \omega^{-1} \, \Tr_{B_1}(L_1, t_1^{+-})
\end{align*}
while
\begin{align*}
\Tr_{B_2}(L_2, t_2^{+-}) &= \beta \,\Tr_{B_2}(L_2', t_2') =   \omega^{-1} \,\Tr_{B_2}(L_2', t_2') , \\
\Tr_{B_2}(L_2, t_2^{-+}) &= \alpha \,\Tr_{B_2}(L_2', t_2')= - \omega^{-5} \,\Tr_{B_2}(L_2', t_2').
\end{align*}

For each of these states $s_1$ for $K_1$, those components of $K_1$ that are not represented on Figure~\ref{fig:Move1bis} and sit below the two arcs shown have the same contribution $X\in \ZZ_{T_1}$ to $\Tr_{T_1}(K_1, s_1)$, while the components of $K_1$ sitting above the two arcs shown contribute $Y\in \ZZ_{T_1}$.

Therefore, the contributions to $\Tr_S(K,s)$ of the four families of states $t_1$, $s_1$, $t_2$  on the left of Figure~\ref{fig:Move1bis} add up to 
\begin{align*}
&\Tr_{B_1}(L_1, t_1^{++})X [Z_{11} Z_{12}^{-1}][Z_{11} Z_{12}]Y  \,\Tr_{B_2}(L_2, t_2^{-+})\\
&\quad\quad + \Tr_{B_1}(L_1, t_1^{-+})X [Z_{11}^{-1} Z_{12}^{-1}] [Z_{11} Z_{12}]Y\,\Tr_{B_2}(L_2, t_2^{-+})\\
& \quad \quad + \Tr_{B_1}(L_1, t_1^{++}) X[Z_{11} Z_{12}][Z_{11} Z_{12}^{-1}]Y \,\Tr_{B_2}(L_2, t_2^{+-})\\
& \quad \quad  + \Tr_{B_1}(L_1, t_1^{+-}) X[Z_{11} Z_{12}][Z_{11}^{-1} Z_{12}^{-1}] Y \,\Tr_{B_2}(L_2, t_2^{+-})\\
=\,&\Tr_{B_1}(L_1, t_1^{++})X (\omega Z_{11} Z_{12}^{-1}) ( \omega Z_{12} Z_{11} ) Y  ( - \omega^{-5}) \,\Tr_{B_2}(L_2', t_2') \\
&\quad\quad + \Tr_{B_1}(L_1, t_1^{-+})X( \omega^{-1} Z_{11}^{-1}  Z_{12}^{-1})(\omega Z_{12}  Z_{11}  )Y ( - \omega^{-5}) \,\Tr_{B_2}(L_2', t_2')\\
& \quad \quad + \Tr_{B_1}(L_1, t_1^{++}) X (\omega^{-1} Z_{11} Z_{12} ) (\omega ^{-1}Z_{12}^{-1} Z_{11}  ) Y\omega^{-1} \,\Tr_{B_2}(L_2', t_2')  \\
& \quad \quad  + \Tr_{B_1}(L_1, t_1^{+-}) X (\omega^{-1}  Z_{11}  Z_{12}  ) (\omega  Z_{12}^{-1}  Z_{11}^{-1}   ) Y  \omega^{-1} \,\Tr_{B_2}(L_2', t_2') \\ 
=&-\omega^{-3}\, \Tr_{B_1}(L_1, t_1^{++}) X Z_{11}^2 Y \,\Tr_{B_2}(L_2', t_2') \\
&\quad\quad -\omega^{-5} \,\Tr_{B_1}(L_1, t_1^{-+})XY  \,\Tr_{B_2}(L_2', t_2')\\
& \quad \quad +\omega^{-3}\, \Tr_{B_1}(L_1, t_1^{++}) X Z_{11}^2 Y \,\Tr_{B_2}(L_2', t_2')  \\
& \quad \quad  + \omega^{-1}\, \Tr_{B_1}(L_1, t_1^{+-}) X Y  \,\Tr_{B_2}(L_2', t_2') \\ 
=&\,
( -\omega^{-5} \,\Tr_{B_1}(L_1, t_1^{-+})  + \omega^{-1}\, \Tr_{B_1}(L_1, t_1^{+-}) )X Y  \,\Tr_{B_2}(L_2', t_2') \\ 
=&\, \Tr_{B_1}(L_1', t_1') X  Y \,\Tr_{B_2}(L_2', t_2')
\end{align*}
which, as required, is the contribution of the right hand side of Figure~\ref{fig:Move1bis} to $\Tr_S(K', s)$. Recall that square brackets $[\enspace ]$ denote here the Weyl quantum ordering (see the end of \S \ref{subsect:MainThm}). Also, note that  the order in which we multiply the contributions of the various components of $K_1$ is really crucial in the above computations.

\begin{figure}[htbp]
\SetLabels
\E(.5 *.8 ) $= $ \\
\E(.5 *.88 ) $? $ \\
\E\L(.9 *.8 ) with $(\epsilon_1, \epsilon_3) \not= (+, -)$ \\
(.245 * .91) $\epsilon_1$\\
(.335 * .53) $\epsilon_3$\\
(.36 * .87) $-$\\
(.4 * .71) $+$\\
(.63 * .91) $\epsilon_1$\\
(.72 * .53) $\epsilon_3$\\
\E(.33 *.2 ) $+ $ \\
\E(.68 *.2 ) $= $ \\
\L\E( 1.05*.2 ) $= \quad0$ \\
\E(.68 *.28 ) $? $ \\
\E(.06 *.28 ) $+ $ \\
\E(.15 *-.04 ) $- $ \\
\E(.18 *.31 ) $+ $ \\
\E(.22 * .14 ) $- $ \\
\E(.44 *.28 ) $+ $ \\
\E(.53 *-.04 ) $- $ \\
\E(.56 *.3 ) $- $ \\
\E(.60 * .14 ) $+$ \\
\E(.81 *.28 ) $+ $ \\
\E(.90 *-.04 ) $- $ \\
\endSetLabels
\centerline{\AffixLabels{ \includegraphics{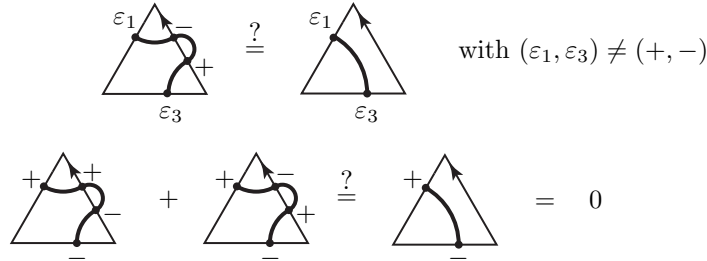} }\hskip 3cm}

\caption{States for Move (II)}
\label{fig:Move2bis}
\end{figure}

We now consider Move~(II). If we again group terms according to compatible states $s_j'$, $t_i'$ for $K'$,  Figure~\ref{fig:Move2bis} lists the possible restrictions to the part of $K$ involved in Move~(II) of all compatible states $s_j$, $t_i$ for $K$ that make non-trivial contributions to $\Tr_S(K,s)$.

In the case of the first line of Figure~\ref{fig:Move2bis}, note that $\Tr_{B_2}(L_2, t_2) = \beta \, \Tr_{B_2}(L_2', t_2')$ by the State Sum Property of Proposition~\ref{prop:BiangleSkein}. As before, let  $X\in \ZZ_{T_1}$ be the contribution of those components of $K_1$ that are not represented on the figure and sit below the two arcs shown, while $Y$ represents the contribution of the components that sits above these two arcs. Then if, as usual, we identify the sign $\epsilon = \pm$ to the number $\epsilon=\pm1$, 
\begin{align*}
\Tr_{T_1} (K_1, s_1) \,\Tr_{B_2}(L_2, t_2)
&= X [Z_{12}Z_{13}^{\epsilon_3}][Z_{11}^{\epsilon_1}Z_{12}^{-1}] Y \beta\, \Tr_{B_2}(L_2', t_2')\\
& = X(\omega^{\epsilon_3} Z_{13}^{\epsilon_3}Z_{12} )(\omega^{-\epsilon_1} Z_{12}^{-1} Z_{11}^{\epsilon_1}) Y \omega^{-1} \,\Tr_{B_2}(L_2', t_2')\\
& = X(\omega^{\epsilon_3 - \epsilon_1-1} Z_{13}^{\epsilon_3} Z_{11}^{\epsilon_1})Y\, \Tr_{B_2}(L_2', t_2')\\
& = X [ Z_{13}^{\epsilon_3} Z_{11}^{\epsilon_1} ] Y\, \Tr_{B_2}(L_2', t_2')\\
&= \Tr_{T_1} (K_1', s_1') \,\Tr_{B_2}(L_2', t_2')
\end{align*}
because $\epsilon_3 - \epsilon_1-1 = -\epsilon_1\epsilon_3$ when $(\epsilon_1, \epsilon_3) \not= (+1, -1)$.

For the second line of Figure~\ref{fig:Move2bis}, there are two families of compatible states $s_1^{-+}$, $t_2^{-+}$ and $s_1^{+-}$, $t_2^{+-}$ for $K$ that correspond to compatible states $s_1'$, $t_2'$ for $K'$ and give non-trivial contributions to $\Tr_{S} (K, s) $. Note that  $s_1'$, $t_2'$  contribute $0$ to $\Tr_{S} (K', s) $. Then,
\begin{align*}
&\Tr_{T_1} (K_1, s_1^{-+}) \,\Tr_{B_2}(L_2, t_2^{-+}) + \Tr_{T_1} (K_1, s_1^{+-}) \,\Tr_{B_2}(L_2, t_2^{+-})\\
& \quad \quad\quad \quad\quad = X [Z_{12}^{-1} Z_{13}^{-1}][Z_{11}Z_{12}] Y \alpha \, \Tr_{B_2}(L_2', t_2')\\
& \quad \quad \quad\quad \quad\quad \quad\quad \quad \quad \quad \quad + X [Z_{12}Z_{13}^{-1}][Z_{11}Z_{12}^{-1}] Y \beta\, \Tr_{B_2}(L_2', t_2')\\
& \quad \quad\quad \quad\quad = X(\omega Z_{13}^{-1}Z_{12}^{-1} )(\omega Z_{12} Z_{11}) Y (-\omega^{-5}) \Tr_{B_2}(L_2', t_2')\\
& \quad \quad\quad \quad\quad \quad \quad \quad \quad + X(\omega^{-1} Z_{13}^{-1}Z_{12})(\omega^{-1} Z_{12}^{-1} Z_{11}) Y \omega^{-1}\, \Tr_{B_2}(L_2', t_2')\\
& \quad \quad\quad \quad\quad = -\omega^{-3}X Z_{13}^{-1}Z_{11} Y  \Tr_{B_2}(L_2', t_2') + \omega^{-3}X Z_{13}^{-1}Z_{11} Y  \Tr_{B_2}(L_2', t_2')\\
& \quad \quad \quad \quad \quad = 0      = \Tr_{T_1} (K_1', s_1') \,\Tr_{B_2}(L_2', t_2').
\end{align*}

This concludes the proof of invariance under Move~(II). 

\begin{figure}[htbp]

\SetLabels
\L\E( .17*.87 ) $= $ \\
\L\E( .38 * .87) $ +$ \\
\L\E( .58*.87 ) $+A^{2} $ \\
\L\E( .78 * .87) $+A^{-2} $ \\
\L\E( .17*.57 ) $= $ \\
\L\E( .38 * .57) $ +$ \\
\L\E( .58*.57 ) $+A^{2} $ \\
\L\E( .78 * .57) $+A^{-2} $ \\
\L\E( .17*.32 ) $= $ \\
\L\E( .38 * .32) $+\quad ( -A^2-A^{-2}+A^{2}+A^{-2})$ \\
\L\E( .17*.08 ) $= $ \\
\endSetLabels
\centerline{\AffixLabels{ \includegraphics{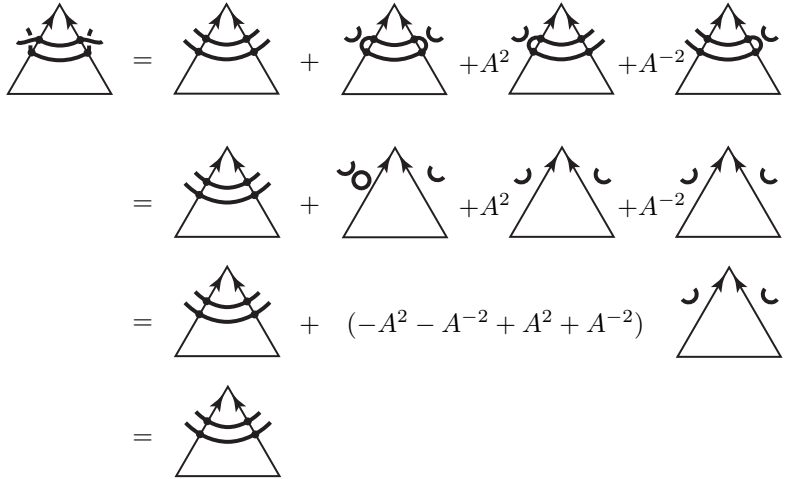} }}

\caption{Skein manipulations for Move~(III)}
\label{fig:Move3bis}
\end{figure}

For Move~(III), instead of state sums, it is more convenient to use the compatibility of the $\Tr_{B_i}(L_i, t_i)$ with the skein relations, as proved by Proposition~\ref{prop:BiangleSkein}. 

The proof of the invariance under Move~(III) is described by Figure~\ref{fig:Move3bis}, where the equalities between linear combinations are understood to apply to the images of the corresponding links under $\Tr_S$. The first equality comes from the fact that $\Tr_{B_1}$ and $\Tr_{B_2}$ are compatible with the skein relations, as proved by Proposition~\ref{prop:BiangleSkein}. The second equality is a consequence of the invariance of $\Tr_S$ under Move~(I), which we just proved. The third equality is a consequence of the fact that, in the skein algebra $\SSS(B_1)$, adding a small unknotted unlinked loop to a skein $[L]$ just multiplies $[L]$ by the scalar $-(A^2 + A^{-2})$ (see Lemma~\ref{lem:UnknotKink}).

\begin{figure}[htbp]

\SetLabels
(.31 * .93) $\epsilon$\\
(.32 * .78) $+$\\
(.42 * .93) $+$\\
(.41 * .78) $+$\\
(.57 * .93) $\epsilon$\\
(.58 * .78) $+$\\
(.68 * .93) $+$\\
(.67 * .78) $+$\\
\E(.5 *.88 ) $= $ \\
\E(.5 *.92 ) $? $ \\
(.31 * .65) $-$\\
(.32 * .53) $-$\\
(.42 * .66) $\epsilon $\\
(.41 * .53) $-$\\
(.57 * .65) $-$\\
(.58 * .53) $-$\\
(.68 * .66) $\epsilon $\\
(.67 * .53) $-$\\
\E(.5 *.62 ) $= $ \\
\E(.5 *.66 ) $? $ \\
(.31 * .4) $\epsilon_1$\\
(.32 * .26) $+$\\
(.43 * .4) $\epsilon_2$\\
(.41 * .27) $-$\\
(.57 * .4) $\epsilon_1$\\
(.58 * .27) $-$\\
(.69 * .4) $\epsilon_2$\\
(.67 * .26) $+$\\
\L(.75 *.35 ) with $(\epsilon_1, \epsilon_2) \not= (-,+)$ \\
\E(.5 *.35 ) $= $ \\
\E(.5 *.39 ) $? $ \\
(.31 * .13) $-$\\
(.32 * -.0) $+$\\
(.42 * .14) $+$\\
(.41 * .005) $-$\\
(.57 * .13) $-$\\
(.58 * .005) $-$\\
(.68 * .14) $+$\\
(.67 * -.0) $+$\\
(.04 * .13) $-$\\
(.06 * .005) $-$\\
(.16 * .14) $+$\\
(.15 * -.00) $+$\\
(.83 * .13) $-$\\
(.85 * -.00) $+$\\
(.95 * .14) $+$\\
(.94 * .005) $-$\\
\E(.5 *.09 ) $= $ \\
\E(.5 *.13 ) $? $ \\
\E(.24 *.09 ) $+ $ \\
\E(.76 *.09 ) $+ $ \\
\endSetLabels
\centerline{\AffixLabels{ \includegraphics{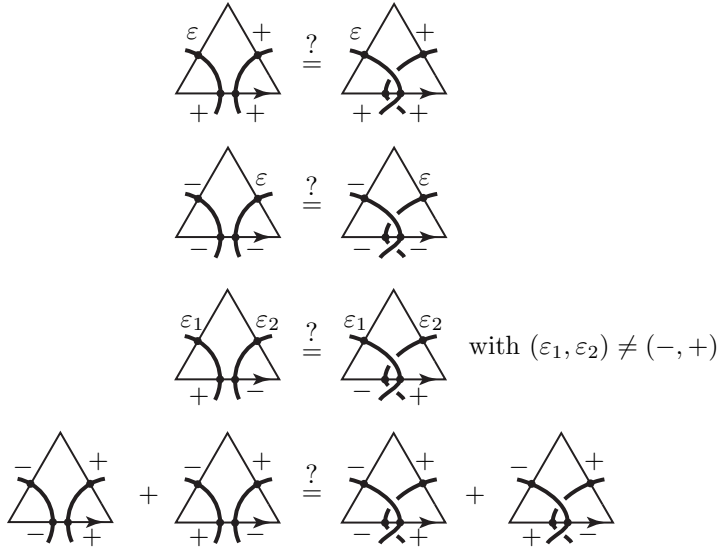} }}

\caption{State sums for Move~(IV)}
\label{fig:Move4bis}
\end{figure}

We now turn to Move (IV). All the states that have a non-trivial contributions are listed in Figure~\ref{fig:Move4bis}. Grouping the contributions of compatible states according to their restrictions outside of the pictures, and according to powers of the generators $Z_{11}$, $Z_{12}$, $Z_{13}$ of $\ZZ_{T_1}$, we have to show the equalities of contributions indicated. 

As usual, let $t_3^{\epsilon\epsilon'}$ be the state for $L_3$ or $L_3'$ where the two boundary points represented are respectively labelled by $\epsilon$, $\epsilon' \in \{ +, -\}$, in this order for the orientation of the edge $T_1\cap B_3$ specified by the arrow (while the value of $t_3^{\epsilon\epsilon'}$ on the other points of $\partial L_3$ and $\partial L_3'$ is determined by the group of compatible states that we are considering).  By  combining Proposition~\ref{prop:BiangleSkein} and Lemma~\ref{lem:RightTwist},
\begin{align*}
\Tr_{B_3}(L_3', t_3^{++}) &= \omega^{-2} \,\Tr_{B_3}(L_3, t_3^{++}) ,\\
\Tr_{B_3}(L_3', t_3^{--}) &= \omega^{-2}  \,\Tr_{B_3}(L_3, t_3^{--}) ,\\
\Tr_{B_3}(L_3', t_3^{+-}) &= \omega^2  \,\Tr_{B_3}(L_3, t_3^{-+}) +(\omega^{-2} - \omega^{6}) \,\Tr_{B_3}(L_3, t_3^{+-}) ,\\
\Tr_{B_3}(L_3', t_3^{-+}) &= \omega^2  \,\Tr_{B_3}(L_3, t_3^{+-}). \\
\end{align*}

Then, for the first line of Figure~\ref{fig:Move4bis}, 
\begin{align*}
\Tr_{T_1}(K_1', s_1')\, \Tr_{B_3}(L_3', t_3')
&=X [Z_{12}Z_{13}]  [Z_{11}^\epsilon Z_{13}] Y \, \Tr_{B_3}(L_3', t_3^{++})\\
&=X \omega^2 [Z_{11}^\epsilon Z_{13}]  [Z_{12}Z_{13}] Y\omega^{-2} \, \Tr_{B_3}(L_3, t_3^{++})\\
&= X [Z_{11}^\epsilon Z_{13}]  [Z_{12}Z_{13}]Y \, \Tr_{B_3}(L_3, t_3^{++})\\
&= \Tr_{T_1}(K_1, s_1)\, \Tr_{B_3}(L_3, t_3),
\end{align*}
where, as usual, $X$ and $Y$ denote the contributions of the components of $K_1$ and $K_1'$ that respectively sit below and above the two arcs represented. 

The case of the second line is almost identical:
\begin{align*}
\Tr_{T_1}(K_1', s_1')\, \Tr_{B_3}(L_3', t_3')
&=X [Z_{12}^\epsilon Z_{13}^{-1}]  [Z_{11}^{-1} Z_{13}^{-1}] Y \, \Tr_{B_3}(L_3', t_3^{--})\\
&=X \omega^2  [Z_{11}^{-1} Z_{13}^{-1}]  [Z_{12}^\epsilon Z_{13}^{-1}]  Y\omega^{-2} \, \Tr_{B_3}(L_3, t_3^{--})\\
&= X [Z_{11}^{-1} Z_{13}^{-1}]  [Z_{12}^\epsilon Z_{13}^{-1}] Y \, \Tr_{B_3}(L_3, t_3^{--})\\
&= \Tr_{T_1}(K_1, s_1)\, \Tr_{B_3}(L_3, t_3).
\end{align*}

For the third  line of Figure~\ref{fig:Move4bis}, 
\begin{align*}
\Tr_{T_1}(K_1', s_1')\, \Tr_{B_3}(L_3', t_3')
&=X [Z_{12}^{\epsilon_2} Z_{13}^{-1}]  [Z_{11}^{\epsilon_1} Z_{13}] Y \, \Tr_{B_3}(L_3', t_3^{-+})\\
&=X \omega^{2(-\epsilon_1\epsilon_2 + \epsilon_2 - \epsilon_1)}   [Z_{11}^{\epsilon_1} Z_{13}]  [Z_{12}^{\epsilon_2} Z_{13}^{-1}]   Y\omega^{2} \, \Tr_{B_3}(L_3, t_3^{+-})\\
&= X [Z_{11}^{\epsilon_1} Z_{13}]  [Z_{12}^{\epsilon_2} Z_{13}^{-1}]   Y \, \Tr_{B_3}(L_3, t_3^{+-})\\
&= \Tr_{T_1}(K_1, s_1)\, \Tr_{B_3}(L_3, t_3)
\end{align*}
as required. Note that $-\epsilon_1\epsilon_2 + \epsilon_2 - \epsilon_1=-1$ exactly when $(\epsilon_1, \epsilon_2) \not= (-1,+1)$. 

The case of the fourth line of  Figure~\ref{fig:Move4bis} is more elaborate. 
\begin{align*}
&\Tr_{T_1}(K_1', s_1^{-+})\, \Tr_{B_3}(L_3', t_3^{-+}) 
 + \Tr_{T_1}(K_1', s_1^{+-})\, \Tr_{B_3}(L_3', t_3^{+-}) \\
& =X [Z_{12} Z_{13}^{-1}] [Z_{11}^{-1} Z_{13}] Y \,\Tr_{B_3}(L_3', t_3^{-+})\\
&\quad \quad \quad\quad \quad \quad\quad \quad \quad + X [Z_{12} Z_{13}] [Z_{11}^{-1} Z_{13}^{-1} ] Y \,\Tr_{B_3}(L_3', t_3^{+-})\\
& = X (\omega Z_{12}Z_{13}^{-1} ) (\omega Z_{13} Z_{11}^{-1}) Y \omega^2 \,\Tr_{B_3}(L_3, t_3^{+-})\\
& \quad \quad \quad \quad + X (\omega^{-1} Z_{12}Z_{13} ) (\omega^{-1} Z_{13}^{-1} Z_{11}^{-1}) Y  \omega^2 \,\Tr_{B_3}(L_3, t_3^{-+})\\
& \quad \quad \quad \quad + X (\omega^{-1} Z_{12}Z_{13} ) (\omega^{-1} Z_{13}^{-1} Z_{11}^{-1}) Y  (\omega^{-2} - \omega^{6})  \,\Tr_{B_3}(L_3, t_3^{+-})\\
& = \omega^4 X Z_{12}  Z_{11}^{-1} Y  \,\Tr_{B_3}(L_3, t_3^{+-}) + X  Z_{12} Z_{11}^{-1} Y \,\Tr_{B_3}(L_3, t_3^{-+})\\
& \quad \quad \quad + \omega ^{-4} X  Z_{12} Z_{11}^{-1} Y  \,\Tr_{B_3}(L_3, t_3^{+-})  - \omega ^4 X Z_{12}Z_{11}^{-1} Y  \,\Tr_{B_3}(L_3, t_3^{+-})\\
&=  \omega^2 X   Z_{11}^{-1}  Z_{12}Y \,\Tr_{B_3}(L_3, t_3^{-+}) + \omega ^{-2} X  Z_{11}^{-1} Z_{12}  Y  \,\Tr_{B_3}(L_3, t_3^{+-})\\
&= X (\omega Z_{11}^{-1}Z_{13}^{-1}) (  \omega Z_{13} Z_{12}) Y \,\Tr_{B_3}(L_3, t_3^{-+})\\
&  \quad \quad \quad \quad \quad \quad \quad \quad \quad + X ( \omega^{-1} Z_{11}^{-1}Z_{13}) ( \omega^{-1} Z_{13}^{-1}Z_{12}) Y \,\Tr_{B_3}(L_3, t_3^{+-})\\
&= X [Z_{11}^{-1}Z_{13}^{-1}] [Z_{12}Z_{13}] Y \,\Tr_{B_3}(L_3, t_3^{-+})\\
&  \quad \quad \quad \quad\quad \quad \quad \quad \quad + X [Z_{11}^{-1}Z_{13}] [Z_{12}Z_{13}^{-1}] Y \,\Tr_{B_3}(L_3, t_3^{+-}).
\end{align*}

This concludes our proof that $\Tr_S(K, s)$ remains invariant under Mover~(IV).

The case of Move~(V) is much simpler. Indeed, by Lemma~\ref{lem:UnknotKink}, 
\begin{align*}
\Tr_{B_1}(L_1', t_1) &= -A^{-3} \,\Tr_{B_1}(L_1, t_1) \\
\Tr_{B_2}(L_2', t_2) &= -A^{3}\, \Tr_{B_2}(L_2, t_2).
\end{align*}
Therefore, when computing $\Tr_S(K', s)$, the two scalars $-A^{-3}$ and $-A^3$ cancel out, and $\Tr_S(K',s) = \Tr_S(K,s)$. 

This concludes our proof  that $\Tr_S(K,s)$ is invariant under the moves (I)--(V),  at least under our original assumption that the biangles $B_1$, $B_2$, $B_3$ touching the triangle $T_1$ where each move takes place are distinct. As indicated at the beginning, we are leaving as an exercise to the reader the task of adapting our arguments to the case where two of these three biangles are equal.

By  Lemma~\ref{lem:GoodPositionMoves}, this concludes the proof of Proposition~\ref{prop:QuantumTraceWellDefined}. 
\end{proof}

\begin{lem}
\label{lem:QuantumTraceDependsSkein}
The above element 
$$
\Tr_S (K, s)  \in \ZZ_\lambda(K),
$$
depends only on the class $[K,s]\in \SSS(S) $ of the framed link $K$ and its state $s$ in the skein algebra. 
\end{lem}

\begin{proof}
We have to show that $\Tr_S$ is compatible with the skein relations, namely that 
$$
\Tr_S(K_1, s) = A^{-1} \Tr_S(K_0, s) + A \Tr_S(K_\infty, s)
$$
when the framed links $K_1$, $K_0$ and $K_\infty$ form a Kauffman triple, namely are related as in Figure~\ref{fig:skein}. 

When we put $K_1$ in good position with respect to the split ideal triangulation $\widehat\lambda$ as in Lemma~\ref{lem:GoodPosition}, we can always arrange that the little ball where $K_1$, $K_0$ and $K_\infty$ differ is located above a biangle $B_j$. For that biangle, Proposition~\ref{prop:BiangleSkein} asserts that $\Tr_{B_j}$ is compatible with the skein relations. In particular, if  $L_j^0$, $L_j^1$, $L_j^\infty$ are the respective intersections of $K_1$, $K_0$ and $K_\infty$ with $B_j\times[0,1]$, it follows from Proposition~\ref{prop:BiangleSkein} that 
$$
\Tr_{B_j}(L_j^1, s_j) = A^{-1}\, \Tr_{B_j}(L_j^0, s_j)  + A\, \Tr_{B_j}(L_j^\infty, s_j) 
$$
for every state $s_j$. By definition of $\Tr_S$ as a state sum, this immediately proves the desired result.
\end{proof}

We are now ready to prove Theorem~\ref{thm:MainThm}. Indeed, the combination of Proposition~\ref{prop:QuantumTraceWellDefined} and Lemma~\ref{lem:QuantumTraceDependsSkein} provides a linear map
$$
\Tr_S \col \SSSS(S)  \to \ZZ_\lambda
$$
defined by $\Tr_S([K,s]) = \Tr_S(K,s)$. This linear map is well-behaved under the superposition operation, so that it is actually an algebra homomorphism. 

Because of its construction as a state sum, it is also immediate that the family of homomorphisms $\Tr_S$ satisfy the State Sum Condition~(1) of Theorem~\ref{thm:MainThm}. 

This State Sum Condition also shows that the homomorphisms $\Tr_S$ are uniquely determined by their restriction to the case where $S$ is a triangle. When $S$ is a triangle, the skein algebra $\SSS(S) $ is generated by simple arcs of the type appearing in Condition~(2) of Theorem~\ref{thm:MainThm} (use the skein relations to eliminate all crossings, and apply Lemma~\ref{lem:UnknotKink} to remove all simple closed curves). The uniqueness part of Theorem~\ref{thm:MainThm} immediately follows.

This concludes the proof of Theorem~\ref{thm:MainThm}. \qed

\section{Invariance under changes of ideal triangulations}

The homomorphism $
\Tr_S \col \SSSS(S)  \to \ZZ_\lambda.
$
provided by Theorem~\ref{thm:MainThm} depends of course on the ideal triangulation $\lambda$ of $S$ considered. We now show that it is well behaved under change of ideal triangulation.

Since we now have to worry about different ideal triangulations, we will write $\Tr_S^\lambda([K,s])$ for the element that we have so far denoted $\Tr_S([K,s]) \in \ZZ_\lambda$

Given two ideal triangulations $\lambda$ and $\lambda'$ of $S$, let $
\Theta_{\lambda\lambda'}^\omega \col \ZZZ_{\lambda'} \to \ZZZ_\lambda
$
be the coordinate change map provided by Theorem~\ref{thm:SquareRootCoordChanges}. 

\begin{thm}
\label{thm:CoordinateChartInvariance}
Given two ideal triangulations $\lambda$ and $\lambda'$ of $S$, and a stated skein $[K,s] \in \SSSS(S) $, the coordinate change map
$$
\Theta_{\lambda\lambda'}^\omega \col \ZZZ_{\lambda'} \to \ZZZ_\lambda 
$$
sends the polynomial $\Tr_S^{\lambda'}([K,s])$ to  the polynomial $\Tr_S^{\lambda}([K,s])$.
\end{thm}

Note that, in general, the coordinate change map $\Theta_{\lambda\lambda'}^\omega $ sends a polynomial $P\in \ZZ_{\lambda'}$ to a \emph{rational fraction} in $\ZZZ_\lambda$. It is therefore surprising that the trace polynomials $\Tr_S^{\lambda}([K,s])$ remain polynomial under coordinate change. 

\begin{proof}[Proof of Theorem~\ref{thm:CoordinateChartInvariance}]
By \cite{Har, Pen}, any two ideal triangulations can be connected to each other by a sequence of diagonal exchanges, as in Figure~\ref{fig:DiagEx}. Since it is proved in \cite[Theorem~25]{Hiatt} that $\Theta_{\lambda\lambda''}^\omega  = \Theta_{\lambda\lambda'}^\omega \circ \Theta_{\lambda'\lambda''}^\omega $ for every three ideal triangulations $\lambda$, $\lambda'$ and $\lambda''$, it will be sufficient to restrict attention to the case where $\lambda$ and $\lambda'$ differ only by a diagonal exchange. 

\begin{figure}[htbp]

\SetLabels
( .167 * 1.05) $\lambda_2 $ \\
( -.04 * .5) $\lambda_5 $ \\
(.167  * -.14) $ \lambda_4$ \\
(.37 * .5) $\lambda_3 $ \\
(.08  * .38) $\lambda_1 $ \\
\E( .5*.5 ) $\longrightarrow $ \\
( .12 * .68) $T_1 $ \\
( .22 * .25) $T_2 $ \\
( .833 * 1.05) $\lambda_2' $ \\
( .63 * .5) $\lambda_5' $ \\
(.833  * -.14) $ \lambda_4'$ \\
(1.03 * .5) $\lambda_3' $ \\
(.92  * .38) $\lambda_1' $ \\
( .88 * .68) $T_1' $ \\
( .78 * .25) $T_2' $ \\
\endSetLabels
\centerline{\AffixLabels{\includegraphics{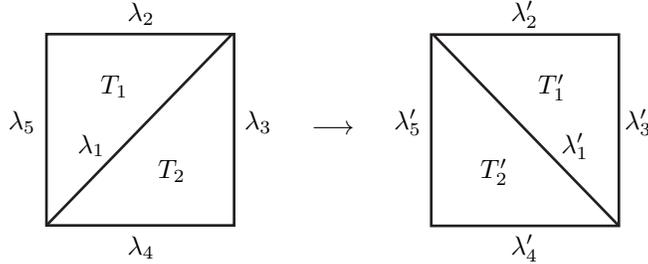}}}\vskip .5cm
\caption{A diagonal exchange}
\label{fig:DiagEx}
\end{figure}

We will assume that the indexing of the edges and faces of $\lambda$ and $\lambda'$ is as in Figure~\ref{fig:DiagEx}. Beware that it is quite possible that there exists identifications between the sides of the square represented, for instance that $\lambda_1=\lambda_2$ or $\lambda_1 = \lambda_3$; however, this will have no impact on our arguments. 

For the split ideal triangulation $\widehat\lambda$ associated to $\lambda$, as usual let $T_j$ be the triangle face associated to the face of $T_j$, and let $B_i$ be the biangle face corresponding to the edge $\lambda_i$ of $\lambda$. We use similar conventions for the split ideal triangulation $\widehat\lambda'$ associated to $\lambda'$. 

Put the framed link $K$ in good position with respect to the split ideal triangulation $\widehat\lambda$, as in Lemma~\ref{lem:GoodPosition}. When doing so, we can always arrange that, above the square $T_1 \cup B_1 \cup T_2$ formed by  the triangles $T_1$, $T_2$ and the biangle $B_1$, the components of $K \cap \bigl((T_1\cup B_1 \cup T_2) \times[0,1]\bigr)$ are all horizontal arcs. Indeed, we can always push any complication of the picture away from the square $T_1\cup B_1\cup T_2$ and into one of the biangles $B_i$ with $i>1$. 

The same property will then hold in $\widehat \lambda'$ since we can always arrange that $T_1' \cup B_1' \cup T_2' = T_1 \cup B_1 \cup T_2$. In particular, $K$ is now in good position with respect to both $\widehat\lambda$ and $\widehat\lambda'$. 

In the state sum expression of $\Tr_S^\lambda ([K,s])$, we can then group the contributions of the components of $K \cap \bigl((T_1\cup B_1 \cup T_2) \times[0,1]\bigr)$ into blocks in $\ZZ_{T_1} \otimes \ZZ_{T_2}$ of one of the following types. 
\begin{enumerate}

\item For components going from $\lambda_2\times[0,1]$ to $\lambda_3 \times [0,1]$:
\begin{enumerate}
\item $[Z_{12} Z_{11}]\otimes [Z_{21} Z_{23}]$;
\item $[Z_{12} Z_{11}]\otimes [Z_{21} Z_{23}^{-1}] \,+\,  [Z_{12} Z_{11}^{-1}]\otimes [Z_{21}^{-1} Z_{23}^{-1} ]$;
\item $[Z_{12} ^{-1} Z_{11}^{-1} ]\otimes [Z_{21} ^{-1} Z_{23}^{-1} ]$.
\end{enumerate}

\item For components going from $\lambda_2\times[0,1]$ to $\lambda_4 \times [0,1]$:
\begin{enumerate}
\item $[Z_{12} Z_{11}]\otimes [Z_{21} Z_{24}] \,+\,  [Z_{12} Z_{11}^{-1}]\otimes [Z_{21}^{-1} Z_{24}] $;
\item $[Z_{12} Z_{11}^{-1}]\otimes [Z_{21}^{-1} Z_{24}^{-1}]$'
\item $[Z_{12}^{-1} Z_{11}^{-1}]\otimes [Z_{21}^{-1} Z_{24}]$;
\item $[Z_{12}^{-1} Z_{11}^{-1}]\otimes [Z_{21}^{-1} Z_{24}^{-1}]$.
\end{enumerate}

\item For components going from $\lambda_2\times[0,1]$ to $\lambda_5 \times [0,1]$:
\begin{enumerate}
\item $[Z_{12} Z_{15}]\otimes 1$;
\item $[Z_{12}^{-1} Z_{15}]\otimes 1$;
\item $[Z_{12}^{-1} Z_{15}^{-1}]\otimes 1$.
\end{enumerate}

\item For components going from $\lambda_3\times[0,1]$ to $\lambda_4 \times [0,1]$:
\begin{enumerate}
\item $1 \otimes [Z_{23} Z_{24}] $;
\item $1 \otimes [Z_{23} Z_{24}^{-1}] $;
\item $1 \otimes [Z_{23}^{-1} Z_{24}^{-1}] $.
\end{enumerate}

\item For components going from $\lambda_3\times[0,1]$ to $\lambda_5 \times [0,1]$:
\begin{enumerate}
\item $[Z_{15} Z_{11}]\otimes [Z_{21} Z_{23}]$;
\item $[Z_{15} Z_{11}]\otimes [Z_{21} Z_{23}^{-1}]$;
\item $[Z_{15}^{-1} Z_{11}]\otimes [Z_{21} Z_{23}]$;
\item $ [Z_{15}^{-1} Z_{11}]\otimes [Z_{21} Z_{23}^{-1}] \,+\,  [Z_{15}^{-1} Z_{11}^{-1}]\otimes [Z_{21} ^{-1}Z_{23}^{-1}] $.
\end{enumerate}

\item For components going from $\lambda_4\times[0,1]$ to $\lambda_5 \times [0,1]$:
\begin{enumerate}
\item $[Z_{15} Z_{11}]\otimes [Z_{21} Z_{24}]$;
\item $[Z_{15}^{-1} Z_{11}]\otimes [Z_{21} Z_{24}] \,+\,  [Z_{15}^{-1} Z_{11}^{-1}]\otimes [Z_{21} ^{-1}Z_{24}]$;
\item $[Z_{15}^{-1} Z_{11}^{-1}]\otimes [Z_{21}^{-1} Z_{24}^{-1}]$.
\end{enumerate}

\end{enumerate}

Then $\Tr_S^{\lambda'}([K,s])$ is obtained from $\Tr_S^{\lambda}([K,s])$ by replacing each of the above blocks by the corresponding block in the list below, while the remaining $Z_{ji}$ with $j>2$ are replaced with the corresponding $Z_{ji}'$. 
\begin{enumerate}

\item For components going from $\lambda_2'\times[0,1]$ to $\lambda_3' \times [0,1]$:
\begin{enumerate}
\item $[Z_{12}'  Z_{13}']\otimes 1$;
\item $[Z_{12}'  Z_{13}' {}^{-1}]\otimes 1$;
\item $[Z_{12}' {}^{-1} Z_{13}' {}^{-1}]\otimes 1$.
\end{enumerate}

\item For components going from $\lambda_2'\times[0,1]$ to $\lambda_4' \times [0,1]$:
\begin{enumerate}
\item $[Z_{12}' Z_{11}']\otimes [Z_{21}' Z_{24}']$;
\item $[Z_{12}' Z_{11}']\otimes [Z_{21}' Z_{24}'{}^{-1} ]$;
\item $[Z_{12}'{}^{-1} Z_{11}']\otimes [Z_{21}' Z_{24}']$;
\item $[Z_{12}'{}^{-1} Z_{11}']\otimes [Z_{21}' Z_{24}' {}^{-1} ] \,+\,  [Z_{12}'{}^{-1} Z_{11}'{}^{-1} ]\otimes [Z_{21}'{}^{-1} Z_{24}'{}^{-1} ]$.
\end{enumerate}

\item For components going from $\lambda_2'\times[0,1]$ to $\lambda_5' \times [0,1]$:
\begin{enumerate}
\item $[Z_{12}' Z_{11}']\otimes [Z_{21}' Z_{25}']$;
\item $[Z_{12}' {}^{-1}Z_{11}']\otimes [Z_{21}' Z_{25}'] \,+\,  [Z_{12}' {}^{-1}Z_{11}'{}^{-1}]\otimes [Z_{21}' {}^{-1}Z_{25}']$;
\item $[Z_{12}'{}^{-1} Z_{11}'{}^{-1}]\otimes [Z_{21}'{}^{-1} Z_{25}'{}^{-1}] $.
\end{enumerate}

\item For components going from $\lambda_3'\times[0,1]$ to $\lambda_4' \times [0,1]$:
\begin{enumerate}
\item $[Z_{12}' Z_{11}']\otimes [Z_{21}' Z_{24}']$;
\item $[Z_{12}' Z_{11}']\otimes [Z_{21}' Z_{24}'{}^{-1}] \,+\,  [Z_{12}' Z_{11}'{}^{-1}]\otimes [Z_{21}' {}^{-1}Z_{24}'{}^{-1}]$;
\item $[Z_{12}' {}^{-1}Z_{11}'{}^{-1}]\otimes [Z_{21}'{}^{-1} Z_{24}'{}^{-1}]$.
\end{enumerate}

\item For components going from $\lambda_3'\times[0,1]$ to $\lambda_5' \times [0,1]$:
\begin{enumerate}
\item $[Z_{13}' Z_{11}' ]\otimes [Z_{21}' Z_{25}' ] \,+\,  [Z_{13}' Z_{11}' {}^{-1}]\otimes [Z_{21}' {}^{-1}Z_{25}' ]$;
\item $[Z_{13}' Z_{11}' {}^{-1}]\otimes [Z_{21}' {}^{-1}Z_{25}' ]$;
\item $[Z_{13}' Z_{11}'{}^{-1} ]\otimes [Z_{21}' {}^{-1}Z_{25}'{}^{-1} ]$;
\item $[Z_{13}' {}^{-1}Z_{11}' {}^{-1}]\otimes [Z_{21}'{}^{-1} Z_{25}'{}^{-1} ]$.
\end{enumerate}

\item For components going from $\lambda_4'\times[0,1]$ to $\lambda_5' \times [0,1]$:
\begin{enumerate}
\item $1 \otimes [Z_{24}' Z_{25}']$;
\item $1 \otimes [Z_{24}' Z_{25}' {}^{-1} ]$;
\item $1 \otimes [Z_{24}'  {}^{-1}Z_{25}'  {}^{-1} ]$.
\end{enumerate}

\end{enumerate}

The coordinate change map $\Theta_{\lambda\lambda'}^\omega \col \ZZZ_{\lambda'} \to \ZZZ_\lambda$ is defined in \cite{Hiatt} by a similar block-by-block analysis. It turns out that it  is very well-behaved with respect to the blocks in the above two lists. Indeed, Hiatt proves in \cite[Lemma~21]{Hiatt} that $\Theta_{\lambda\lambda'}^\omega$ sends each block of $\ZZ_{T_1'}\otimes \ZZ_{T_2'}$ in the list above to the corresponding block  in $\ZZ_{T_1}\otimes \ZZ_{T_2}$ in the first list, while sending each element of $\ZZ_{T_j'}$ with $j>2$ to the element of $\ZZ_{T_j}$ obtained by removing the primes ${}'$. When combined with our original observations, this proves that 
$$
\Theta_{\lambda\lambda'}^\omega \bigl( \Tr_S^{\lambda'}([K,s]) \bigr) = \Tr_S^{\lambda}([K,s])
$$
in the case considered, namely when $\lambda$ and $\lambda'$ differ by a diagonal exchange. 

As observed at the beginning of our discussion, this completes the proof of Theorem~\ref{thm:CoordinateChartInvariance}.
\end{proof}

\section{Injectivity of the quantum trace homomorphism}

We conclude with a simple observation.

\begin{prop}
\label{prop:QTraceInjective}
The quantum trace homomomorphism
$$\Tr_\lambda^\omega  \col \SSS(S) \to \ZZZ_\lambda$$
of Theorem~{\upshape\ref{thm:MainThmIntro}} is injective.
\end{prop}

\begin{proof}
As a vector space, the skein algebra is clearly generated by the family of all skeins $[K]\in \SSS(S)$ that are \emph{simple}, in the sense that they are represented by $1$--submanifolds of $S$ (with no crossings, and with vertical framing) whose components are not homotopic to 0.

For such a simple skein $[K]$, our state sum construction of the Laurent polynomial $\Tr_\lambda^\omega  ([K]) \in  \ZZZ_\lambda$ shows that its highest degree term  is a non-zero scalar multiple of $Z_1^{k_1} Z_2^{k_2}\dots Z_n^{k_n}$, where $k_i\geq 0$ is the geometric intersection number of $K$ with the $i$--th edge $\lambda_i$ of the ideal triangulation~$\lambda$. 

The key observation is now that a simple skein $[K]$ can be completely recovered from the collection $(k_1, k_2, \dots, k_n)$ of its geometric intersection numbers. It easily follows that  the image under $\Tr_\lambda^\omega$ of a non-trivial linear combination of simple skeins cannot be 0 (focus attention on a term for which $(k_1, k_2, \dots, k_n)$ is maximal), which proves that the kernel of $\Tr_\lambda^\omega$ is trivial.  
\end{proof}
Incidentally, the above argument also provides another proof that simple skeins are linearly independent in $\SSS(S)$.

\end{document}